\documentclass[numbers,webpdf,imaiai]{ima-arXiv}

\usepackage{listings}
\usepackage{color}
\usepackage{esint}
\usepackage{cleveref}
\usepackage{float}
\usepackage{amsmath}
\usepackage{caption}
\usepackage{bbm}

\usepackage{hyperref}

\def\R{\mathbb{R}}

\def\newmap{{single-slice matching~}}
\def\Mnewmap{{matrix-slice matching~}}

\def\W{{\mathcal{W}}}
\def\leb{{\mathcal{L}^n}}

\newcommand{\bx}{{{x}}}

\newcommand{\bth}{{{\theta}}}

\newcommand{\id}{\operatorname{id}}
\newcommand{\E}{{\mathbb{E}}}

\newcommand{\Ttheta}{{T^{\theta}}} 
\newcommand{\Tthetaone}{{T^{\theta_1}}} 
\newcommand{\Tthetatwo}{{T^{\theta_2}}}
\newcommand{\Tthetai}{{T^{\theta_i}}} 
\newcommand{\Tthetan}{{T^{\theta_n}}}
\newcommand{\singleslice}{{T_{\sigma,\theta}}}
\newcommand{\singleslicek}{{T_{\sigma,\theta_k}}}
\newcommand{\Pslice}{{T_{\sigma,P}}}
\newcommand{\Pjslice}{{T^j_{\sigma,P}}}
\newcommand{\Pslicek}{{T_{\sigma_k,P_k}}}
\newcommand{\Pkjslice}{{T^j_{\sigma_k,P_k}}}
\newcommand{\Pkplusjslice}{{T^j_{\sigma_{k},P_k}}}
\newcommand{\jkPslice}{{T^j_{\sigma_k,P}}}
\newcommand{\Tkthetaik}{{T_k^{\theta^k_{i}}}}
\newcommand{\Tkplusthetaik}{{T_{k+1}^{\theta^k_{i}}}}

\theoremstyle{thmstyletwo}%
\newtheorem{theorem}{Theorem}
\newtheorem{proposition}[theorem]{Proposition}%

\newtheorem{example}{Example}%
\newtheorem{remark}{Remark}%
\newtheorem{lemma}{Lemma}[section]
\newtheorem{assumption}{Assumption}
\newtheorem{corollary}{Corollary}

\newtheorem{definition}{Definition}

\numberwithin{equation}{section}

\begin{document}

\DOI{DOI HERE}
\copyrightyear{2024}
\vol{00}
\pubyear{2024}
\access{Advance Access Publication Date: Day Month Year}
\appnotes{Paper}
\copyrightstatement{Published by Oxford University Press on behalf of the Institute of Mathematics and its Applications. All rights reserved.}
\firstpage{1}


\title[Measure transfer via stochastic slicing and matching]{Measure transfer via stochastic slicing and matching}

\author{Shiying Li*
\address{\orgdiv{Department of Mathematics}, \orgname{University of Nebraska - Lincoln}, \orgaddress{\street{Avery Hall}, \postcode{68508}, \state{Nebraska}, \country{USA}}}
}

\author{Caroline Moosm\"{u}ller
\address{\orgdiv{Department of Mathematics}, \orgname{University of North Carolina at Chapel Hill}, \orgaddress{\street{ 3250 Phillips Hall}, \postcode{27599}, \state{North Carolina}, \country{USA}}}}

\author{{Yongzhe Wang 
\address{\orgdiv{Department of Mathematics}, \orgname{University of Pennsylvania}, \orgaddress{\street{David Rittenhouse Laboratory}, \postcode{19104}, \state{Pennsylvania}, \country{USA}}}}
}

\authormark{Shiying Li, Caroline Moosm\"{u}ller and Yongzhe Wang}

\corresp[*]{Corresponding author: \href{email:email-id.com}{sli82@unl.edu}}

\abstract{This paper studies iterative schemes for measure transfer and approximation problems, which are defined through a slicing-and-matching procedure. Similar to the sliced Wasserstein distance, these schemes benefit from the availability of closed-form solutions for the one-dimensional optimal transport problem and the associated computational advantages. While such schemes have already been successfully utilized in data science applications, not too many results on their convergence are available. The main contribution of this paper is an almost sure convergence proof for stochastic slicing-and-matching schemes. The proof builds on an interpretation as a stochastic gradient descent scheme on the Wasserstein space. Numerical examples on step-wise image morphing are demonstrated as well.}


\maketitle

\section{Introduction}

Optimal transport and the Wasserstein distance have gained widespread interest in the machine learning community, as they provide a natural framework for dealing with datasets consisting of point clouds or measures. Some areas of application include generative models \cite{arjovsky2017wasserstein}, (semi-supervised) learning \cite{solomon2014wasserstein,frogner2015WassersteinLoss}, signal processing \cite{kolouri2017optimal}, and imaging \cite{rubner2000earth}.

The optimal transport problem seeks to find the best way (in the sense of cost-minimizing) to transport one measure into another \cite{Villani1}. In the Monge formulation, one seeks a map $T$ which transports the original measure $\sigma$ to the target measure $\mu$ (i.e. $T_{\sharp}\sigma = \mu$) and minimizes the cost
\begin{equation}\label{intro:W2}
    W_2^2(\sigma,\mu):=\min_{T: T_{\sharp}\sigma = \mu} \int_{\mathbb{R}^n} \|T(x)-x\|^2 \, d\sigma(x).
\end{equation}
Here $W_2$ is the $2$-Wasserstein distance and the argmin is the optimal transport map from $\sigma$ to $\mu$. Note that in this set-up, a solution to \eqref{intro:W2} may not exist; \cite{brenier1991} provides conditions for existence (such as absolute continuity of $\sigma$), and Kantorovich \cite{kantorovich42} relaxed this framework, seeking joint distributions rather than maps (see \Cref{sec:prelims} for more details). We introduce the Monge setting here, as the main focus of this paper is studying transportation maps, which mimic certain behaviors of the optimal transport map.
We also mention that \eqref{intro:W2} is a special case of $p$-Wasserstein distances, and the Euclidean distance can also be replaced by more general cost functions \cite{Villani1,mccann01}.

While very successful in applications, optimal transport can be computationally expensive, especially in high dimensions. The problem \eqref{intro:W2} becomes a linear program of computational order $O(m^3\log(m))$, {where $m$ is the number of samples in the discrete setting.}
For this reason, there is  
interest in approximation schemes both for the Wasserstein distance as well as for the optimal transport map. A well-known approach is entropic regularized optimal transport (Sinkhorn distances) \cite{cuturi-2013}, which significantly reduces the computational cost by using matrix scaling algorithms \cite{sinkhorn67}. Other approximation schemes take advantage of particular properties of the underlying set of measures; linearized optimal transport, for example, uses linear distances in a tangent space, which approximate the Wasserstein distance if the set of measures has an almost flat structure \cite{park2018cumulative, moosmueller2020linear, wang2013, aldroubi20}. In this paper, we are interested in a particular type of approximation, namely \emph{sliced Wasserstein distances} \cite{bonneel2015sliced,rabin2010ShapeRetrieval,rabin2012wasserstein,mahey2023slicedGeodesics}, which make use of the fact that Wasserstein distances in one dimension can be computed easily through the cumulative distribution functions (CDFs). More concretely, we are interested in the slicing idea underlying these distances, which can be used to construct transport maps. These maps, in turn, give rise to iterative schemes for measure approximation. We introduce these ideas in the next section.

\subsection{Sliced Wasserstein distances and iterative approximation schemes}

The sliced Wasserstein distance \cite{bonneel2015sliced,rabin2010ShapeRetrieval} between two measures 
$\sigma, \mu$ is given by
\begin{equation*}
    SW_2^2(\sigma,\mu) = \int_{S^{n-1}} W_2^2(\sigma^{\bth},\mu^{\bth}) \, du(\theta),
\end{equation*}
where $\mu^{\bth}$ is the one-dimensional measure defined by $\mu^{\bth} ={\mathcal{P}_{{\theta}}}_{\sharp}\mu$ with ${\mathcal{P}_{{\theta}}}(x)= \left<x,\theta\right> = x\cdot \theta$ the projection onto the unit vector $\theta$. Here $u$ denotes the uniform probability measure over $S^{n-1}$. The Wasserstein distance $W_2$ under the integral is between one-dimensional measures, and can therefore be computed by the CDFs of $\sigma^{\bth}$ and $\mu^{\bth}$ (no optimization necessary). {Under certain regularity assumptions on a measure $\eta$, the ``sliced" measure $\eta^{\theta}$  can be viewed as the Radon transform of $\eta$ evaluated at $\theta$, see for example \cite[Proposition 6]{bonneel2015sliced}.}

Due to its simplicity, the sliced Wasserstein distance is often used as a replacement for the full Wasserstein distance, and has been successful in many applications, such as texture mixing and barycenters \cite{rabin2012wasserstein}, shape retrieval \cite{rabin2010ShapeRetrieval}, neural style transfer \cite{li2022neural} and radiomics studies \cite{belkhatir2022texture}.

The sliced Wasserstein distance has also been extended to different settings, such as unbalanced and partial transport problems \cite{sejourne2023unbalanced,bonneel2019SPOT,Bai2023partial} and generalized slicing \cite{kolouri2019generalizedSlicesW,Kolouri2022ICASSP}.

Closely related to sliced Wasserstein distances is the idea of using slices to define transport maps. In the simplest setting {where $\mu$ is fixed} \cite{pitie2007automated}, one chooses a line $\theta \in S^{n-1}$, and defines
\begin{equation}\label{intro:single-Sliceform}
 \singleslice(\bx) := \bx + (\Ttheta(x\cdot \theta)-x\cdot \theta)\,\theta,
\end{equation}
where $\Ttheta$ is the one-dimensional optimal transport map between the slices $\sigma^{\theta},\mu^{\theta}$, which can again be explicitly computed through the CDFs.
Similarly, one can use multiple directions $\theta_i$ to define a map in the spirit of \eqref{intro:single-Sliceform}. For example, \cite{pitie2007automated} uses an orthogonal matrix $P= [\theta_1,\ldots,\theta_n]$ and defines
\begin{equation}\label{intro:Sliceform}
   \Pslice(\bx) := \bx + P\begin{bmatrix}\Tthetaone(\bx\cdot \bth_1)-\bx\cdot \bth_1\\\Tthetatwo(\bx\cdot \bth_2)-\bx\cdot \bth_2\\ \vdots\\\Tthetan(\bx\cdot \bth_n)-\bx\cdot \bth_n
    \end{bmatrix}.
\end{equation}

While \eqref{intro:single-Sliceform} and \eqref{intro:Sliceform} are not necessarily the optimal transport map between $\sigma$ and $\mu$, they can be used in approximation schemes. Similar to efforts in \cite{baptista2023ApproxFramework}, the Knothe-Rosenblatt rearrangement \cite{rosenblatt52arrangement,Knothe1957convex}, and no-collision transport maps \cite{negrini2024applications,nurbekyan2020nocollision}, these maps can be used as an easy-to-compute replacement for the actual optimal transport map or to approximate $\mu$ itself.
{Both Knothe-Rosenblatt (KR) and no-collision aim to produce a transport (mass-preserving) map between the measures. 
No-collision builds an injective map by recursively splitting the domain with hyperplanes and, in each cell, matching the one-dimensional distributions along the split direction (that is, matching marginals conditioned on the side of the hyperplane).
Slice-matching does not partition the domain: it projects onto prescribed directions and applies one-dimensional transports along those projections, aligning the projected marginals in parallel, but a single pass typically does not define a full transport to the target, so iteration is used.
}
{On the other hand,  when the slicing directions are canonical, e.g., $P=I$, the slice-matching map reduces to a diagonal map, which  a special case of the triangular structure used by KR construction. In the KR construction, $T_m$ transports the conditional distribution of $\mu$ given $x_{1:i-1}$ to the conditional distribution of $\mu$ given $y_{1:i-1}$, with $y_{1:i-1}=T_{1:i-1}(x_{1:i-1})$. In contrast, slice-matching works in arbitrary orthonormal directions: each $f_i$ transports the marginal distribution of $\sigma$ along the direction $\theta_i$ to the corresponding marginal of $\mu$ along $\theta_i$. In the canonical basis ($P=I$), the diagonal slice map and the KR map coincide if and only if, for each $i$, $\sigma^{i\mid x_{1:i-1}}=\sigma^{\theta_i}$ and $\mu^{i\mid y_{1:i-1}}=\mu^{\theta_i}$ ($\sigma^{i\mid x_{1:i-1}} \textrm{~and~} \mu^{i\mid y_{1:i-1}}$ are the conditional measures), i.e., when both measures have independent coordinates in that basis.}

\begin{figure}
    \centering
    \includegraphics[width=0.7\linewidth]{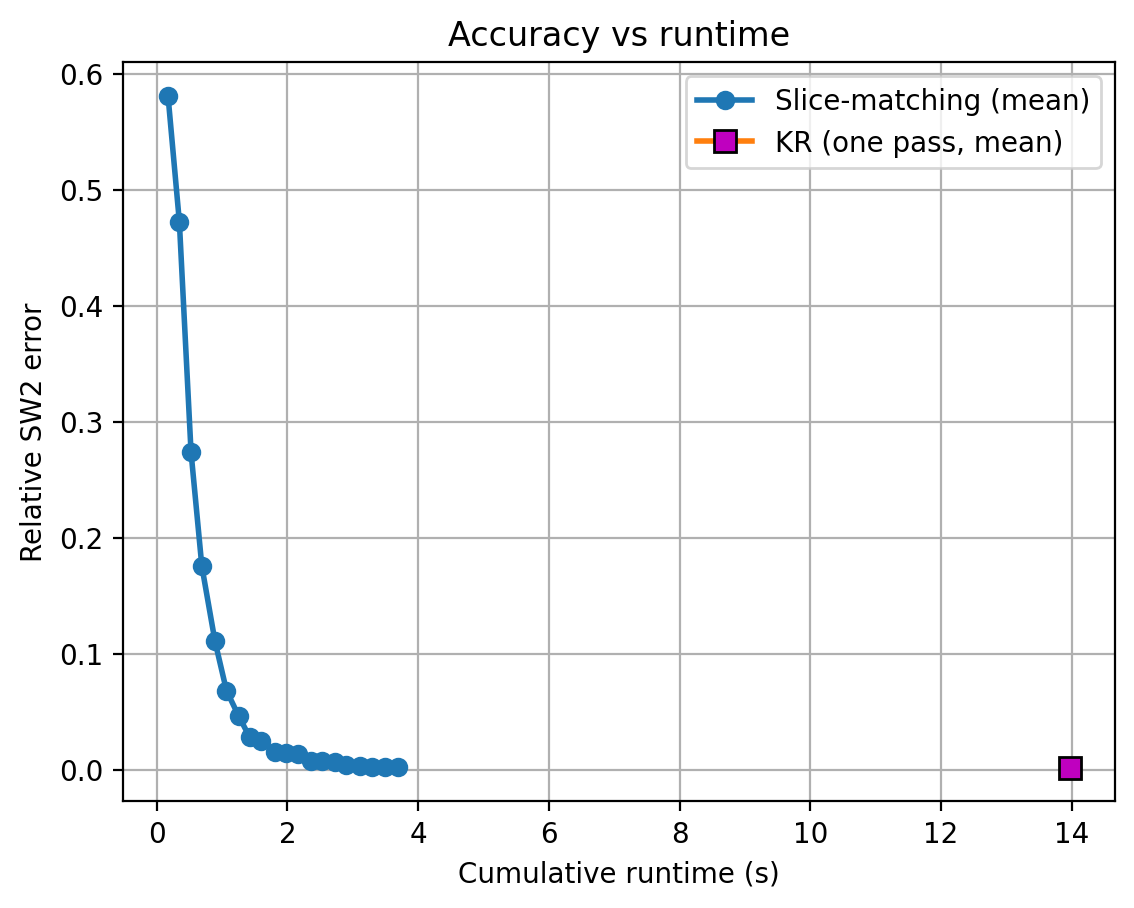}
    \caption{Comparison of measure transfer using slice-matching iterative scheme and Knothes-Rosenblatt transport (Experiment \ref{Exp:rotGaussian}.1). Source: $\mathcal{N}(0, \Sigma)$, where $\Sigma = \textrm{diag}(4, 1, 1/4)$. Target: $\mathcal{N}(0, Q\Sigma Q^T)$, where $Q$ is a random orthogonal matrix. Each distribution is approximated using  $N = 200{,}000$ i.i.d. samples. The plot shows the relative sliced-Wasserstein error (averaged over 5 trials) versus  cumulative time.}
    \label{KRfig:n3rotatedGau-relError}
\end{figure}

{
In high dimensions, slice-matching can be attractive from a computational perspective: each map \(f_k\) acts only on a one-dimensional marginal along direction \(\theta_k\), allowing these maps to be computed independently and in parallel. The trade-off is that one application generally does not transport \(\sigma\) to \(\mu\), so multiple iterations are used for approximation, which directly motivates the analysis carried out in this manuscript. By contrast, the KR map produces a transport in a single step, but each component \(T_i : \mathbb{R}^i \to \mathbb{R}^i\) depends on the previously computed coordinates through \(y_{1:i-1} = T_{1:i-1}(x_{1:i-1})\). In practice, evaluating \(T_k\) requires computing conditional CDFs \(F_\sigma(\cdot \mid x_{1:i-1})\) and \(F_\sigma(\cdot \mid y_{1:i-1})\) for many different conditioning values \((x_{1:i-1}, y_{1:i-1})\), and in a sequential manner. Unless the densities have a special factorizable form \cite{baptista2023representation}, this conditional evaluation can be considerably more involved than marginal estimation in slice-matching. For a numerical comparison, we plot---averaged over five trials---the relative sliced-Wasserstein error of the transported measure obtained by iterative applications of the slice-matching maps versus the cumulative runtime;  we also plot the corresponding KR error versus its (one-pass) runtime. \Cref{KRfig:n3rotatedGau-relError} reports the experiment where the source is an anisotropic Gaussian distribution in $\R^3$ and the target is a rotated version of the source (cf. \Cref{sec:KRcomparison} for more details). As expected, the KR attains a small error in a single application,  whereas the slice matching requires multiple iterations to reach comparable errors; nonetheless, for a given error,  the cumulative runtime of slice matching is much smaller than that of KR.    }

To this end, \cite{pitie2007automated} suggested to iteratively apply maps of the form \eqref{intro:Sliceform} to $\sigma$, with the aim of approximating $\mu$ in some distance (\cite{pitie2007automated} studies KL-divergence, we are concerned with the Wasserstein and the sliced Wasserstein distance). In particular, \cite{pitie2007automated} proposes the following iterative scheme: Choose a random sequence of orthogonal matrices $\{P_k\} \subset O(n)$ and let $\sigma_0=\sigma$. Then
\begin{equation}\label{intro:iter}
 \sigma_{k+1} = ({\Pslicek})_{\sharp}\sigma_k, \quad k\geq 0.
\end{equation}
The main motivation in \cite{pitie2007automated} is color transfer of images, but this idea generalizes to a myriad of applications, including texture mixing \cite{rabin2012wasserstein}, shape retrieval \cite{rabin2010ShapeRetrieval}, barycenter problems \cite{bonneel2015sliced}, sampling \cite{paulin2020}, generative modeling with normalizing flows \cite{dai2021sliced}, in addition to many interesting connections to gradient descent and Procrustes analysis \cite{zemel2019frechet} as well as  various gradient flows of probability measures \cite{liutkus2019slicedFlows, bonet2022efficient}, see \Cref{Sec:rel_GD}.

To justify the use of iterative scheme of this type in data science applications, convergence results of the form $\sigma_k \to \mu$ as $k\to \infty$ are crucial. Currently, not too many results in this direction are available: \cite{pitie2007automated} shows convergence of \eqref{intro:iter} in KL-divergence, when both $\sigma$ and $\mu$ are Gaussians. These results were further refined by \cite{bonnotte13thesis}.
In addition, \cite{meng2019large} shows convergence of estimators in the Wasserstein distance, i.e.\ when the number of samples goes to infinity, assuming that the number of iterations $k$ is large enough (suitable sizes of $k$ depend on the dimension of the space). 
The present manuscript aims at contributing towards these efforts by providing a rigorous convergence proof for a variant of the scheme \eqref{intro:iter}, which we introduce in the next section.

\subsection{Stochastic iterative approximation schemes}
In this paper, we are interested in a stochastic version of the scheme \eqref{intro:iter}, namely
\begin{equation}\label{intro:randomized_iter}
  \sigma_{k+1} = ((1-\gamma_k)\operatorname{id}+\gamma_k{\Pslicek})_{\sharp}\sigma_k, \quad k\geq 0,
\end{equation}
where $\gamma_k \in [0,1]$ is a sequence satisfying the classical {stochastic} gradient descent assumptions {\cite{robbins1951stochastic}} $ \sum_{k=0}^{\infty} \gamma_k = \infty$ and $\sum_{k=0}^{\infty} \gamma_k^2 < \infty$. Stochasticity is obtained by choosing $P_k$ as i.i.d.\ samples from the Haar measure on $O(n)$.
We mention that the original scheme \eqref{intro:iter} does not fall into this class of iteration schemes, since the constant step-size $\gamma_k=1$ does not satisfy the assumptions. While our convergence results, as outlined in \Cref{sec:main_contributions}, hold for a large class of iteration schemes, they do not hold for \eqref{intro:iter}.

The version \eqref{intro:randomized_iter} was first studied in \cite{bonneel2015sliced,rabin2010ShapeRetrieval,rabin2012wasserstein}, where the main focus was on applications and Wasserstein barycenters. Further significant contributions were made in \cite{bonnotte13thesis}, as will be outlined throughout this manuscript. 
These papers furthermore observe that \eqref{intro:randomized_iter} can be interpreted as a stochastic gradient descent scheme for a loss function closely related to
\begin{equation}\label{intro:loss}
   {H}(\sigma) = \frac12 SW_2^2(\sigma,\mu).
\end{equation}
The main focus of \cite{bonneel2015sliced,rabin2010ShapeRetrieval,rabin2012wasserstein} is point cloud data, which translates the iteration \eqref{intro:randomized_iter} into a stochastic gradient descent scheme on $\mathbb{R}^N$. Our paper is concerned with measures and therefore uses stochastic gradient descent schemes in the Wasserstein space, building on the recent results of {\cite{backhoff2025stochastic}}. We summarize our main contribution in the next section.

\subsection{Main contribution}\label{sec:main_contributions}

The main contribution of this paper is an almost sure (a.s.) convergence proof for the stochastic iterative scheme \eqref{intro:randomized_iter}, when $\sigma$ and $\mu$ are measures rather than point clouds and slices are i.i.d.\ drawn from the Haar measure on $O(n)$. The proof uses stochastic gradient descent on the Wasserstein space with a modified version of the loss \eqref{intro:loss}.
Our result is motivated by the observations in \cite{bonneel2015sliced,rabin2010ShapeRetrieval,rabin2012wasserstein,bonnotte13thesis} on point cloud data and the proof techniques of the recent paper \cite{backhoff2025stochastic}, which studies stochastic gradient descent schemes and population barycenters in the Wasserstein space. 
In particular, we derive the following result:
\begin{theorem}[Special case of \Cref{mainthm-Pj}; informal version]\label{intro:mainthm-Pversion}
    Consider two measures $\sigma_0,\mu$ over $\R^n$ and let $P_k \overset{\textrm{i.i.d}}{\sim} u_n, k\geq 0$,  where $u_n$ is the Haar probability measure on $O(n)$. Define $\sigma_k$ by the iteration \eqref{intro:randomized_iter} using $\sigma_0,\mu$ and $P_k$.
       Then under some technical assumptions we get
    \begin{equation*}
        \sigma_k \xrightarrow{W_2} \mu, \quad  \text{a.s.\ as } k\to \infty.
    \end{equation*} 
\end{theorem}
In fact, we show an a.s.\ convergence result for any ``inbetween'' slicing scheme using $1\leq j\leq n$ slices, and \Cref{intro:mainthm-Pversion} is a special case with $n$ slices. Our result also includes the ``single-slice scheme'', which uses the map $\singleslice$ from \eqref{intro:single-Sliceform} in the iteration with i.i.d.\ samples $\theta_k$ drawn from the uniform probability measure on $S^{n-1}$ (\Cref{{mainthm}}).

The main reason why the proof techniques of \cite{backhoff2025stochastic} are useful for our setting is a reformulation to a special type of barycenter problem, as we outline in the proof of \Cref{mainthm-Pj}. {{
We also note that our convergence proof relies on a strong hypothesis related to the local uniqueness of Karcher means \cite{zemel2019frechet, backhoff2025stochastic} near the minimizer of an associated Fr\'echet functional (the Wasserstein barycenter objective); see Assumption (A1-ii) for more details. At a high level,  this means that the functional on the space of probability measures admits a single local critical point. In our setting, the functional $H$  in \eqref{intro:loss} is  defined using the sliced Wasserstein distance with a fixed target. A general theory about the uniqueness of Karcher means is an interesting and challenging problem on its own and remains largely open, as also noted in \cite{backhoff2025stochastic, zemel2019frechet}, yet it is highly relevant to the analysis of gradient-descent-based methods for functional minimization. Current results establish uniqueness in special settings, such as a finite collection of measures \cite{zemel2019frechet} and  parametric pushforward families of a fixed base measure under strong regularity assumptions \cite{bigot2018characterization}.

} }

\subsection{Relation to gradient flows}\label{Sec:rel_GD}
In the case of functions $F:\R^n \to \R$, a gradient flow equation is of the form
\begin{equation}\label{eq:GF-Rn}
    \dot{x} = - \nabla F (x).
\end{equation}
The implicit Euler scheme for \eqref{eq:GF-Rn}   can be reformulated as a ``minimizing movement scheme''\cite{de1993new}, which can then be generalized to a scheme operating on measures rather than points in Euclidean space (replacing Euclidean distance by, for example, the Wasserstein distance) \cite{benamou2000Fluid,Jordan1998JKO}. This gives the $W_2$-gradient flow scheme \cite{ambrosio2008gradient,santambrogio2017flows}. 

In this paper we are interested in problems of the form $\min_{\sigma \in \mathcal{W}_2(\R^n)}\mathcal{F}(\sigma)$, which, if formulated for $F:\R^n \to \R$, could be tackled by a gradient descent scheme 
\begin{equation}\label{eq:GD-Rn}
    x_{n+1} = x_n -h \nabla F(x_n),
\end{equation}
i.e.\ by the explicit Euler scheme for \eqref{eq:GF-Rn}. For $\mathcal{F}: \mathcal{W}_2(\R^n) \to \R$, we consider the following version of \eqref{eq:GD-Rn}:
\begin{equation*}
    \sigma_{n+1} = \left(\id -h \mathcal{F}^{\prime}\right)_{\sharp}\sigma_n,
\end{equation*}
where $\mathcal{F}^{\prime}$ is a (formal) Fr\'echet derivative of $\mathcal{F}$ \cite{ , ambrosio2008gradient,zemel2019frechet}. This provides an interpretation of gradient descent on the Wasserstein space \cite{zemel2019frechet, Molchanov02_steepest}, and shows the close relation of this idea to gradient flows.

In the case of empirical measure, \cite[Proposition 5.2.7]{bonnotte13thesis} shows that the iterative scheme \eqref{intro:iter} can be interpreted as a kind of gradient descent scheme for the functional 
$\mathcal{F}_{P_k}(\sigma) = \sum_{i=1}^nW_2^2(\sigma^{\theta_i^k}, \mu^{\theta_i^k})$ with orthogonal matrix $P_k = [\theta_1^k,\cdots,\theta_n^k]$ and step-size $h=1$. The functional, however, depends on the iteration variable $k$. To deal with this issue, \cite{rabin2012wasserstein} suggested to integrate over $O(n)$, i.e.\ to consider the functional
\begin{equation}\label{eq:integrate-loss}
    \mathcal{F}(\sigma) = \int_{O(n)}\sum_{i=1}^n W_2^2(\sigma^{\theta_i}, \mu^{\theta_i})\,dP,
\end{equation}
which is closely related to the loss $L$ of \eqref{intro:loss} and hence the sliced Wasserstein distance. \cite{rabin2012wasserstein} furthermore studies the gradient flow related to this functional. Motivated by the recent results in \cite{backhoff2025stochastic}, we give a stochastic gradient descent interpretation for the iterative scheme \eqref{intro:randomized_iter} in terms of the functional \eqref{eq:integrate-loss}. While \eqref{intro:randomized_iter} can be alternatively viewed as a batch gradient descent procedure with respect to the functional \eqref{intro:loss}, see \cite{bonneel2015sliced}, a version of \eqref{eq:integrate-loss} conveniently leads to a unifying stochastic gradient descent convergence analysis for a more general framework of iterative slice-matching schemes.

Other work in this area include: $W_2$-gradient flows for functionals defined using generalized sliced probability metrics \cite{Kolouri2022ICASSP,kolouri2019generalizedSlicesW}, for barycenter problems using functionals of the form {$\int_{\W_2(\R^n)} W_2^2(\sigma,\mu) d\Pi(\mu)$ for some $\Pi$} \cite{,chewi2020BuresWasserstein}, for $SW_2$ with entropy functionals with applications in generative modeling \cite{liutkus2019slicedFlows}, and $SW_2$-gradient flows (replacing the Euclidean distance by the sliced Wasserstein distance) for general functionals \cite{bonet2022efficient}.

\subsection{Structure of the paper}
\Cref{sec:prelims} summarizes the necessary {preliminaries} such as optimal transport, Wasserstein distance and the slicing procedure, which can be used to define the sliced Wasserstein distances as well as ``sliced'' transport maps. In \Cref{sec:slice-matching} we introduce a generalized form of the slice-matching maps of \cite{pitie2007automated} and show some basic properties of these transports, including a connection to \emph{compatible maps} as introduced in \cite{khurana2022supervised,aldroubi20}. \Cref{sec:iterative_scheme} then shows the relation of a generalized version of the iterative scheme \eqref{intro:randomized_iter} to stochastic gradient descent on the Wasserstein space. The main result on a.s.\ convergence of this scheme is stated in \Cref{mainthm-Pj}, along with a corollary showing convergence for the ``single-slice'' scheme, which is based on iterating \eqref{intro:single-Sliceform}.  The proofs are made available in \Cref{appendixProofs}. The paper closes with numerical experiments on morphing images in \Cref{sec:experiments}.

\section{Preliminaries}\label{sec:prelims}

\subsection{Optimal transport preliminaries}
By $\mathcal{P}(\R^n)$ we denote the space of {Radon} probability measures on $\R^n$, and by $\mathcal{P}_{ac}(\R^n)$ the space of absolutely continuous  {Radon} probability measures (with respect to the Lebesgue measure).

We are furthermore interested in the quadratic Wasserstein space, which is the space of probability measures $\sigma$ with finite second moment
$\int_{\R^n}\|x\|^2d\sigma(x)<\infty$. We denote this space by $\W_2(\R^n)$. We also define $\W_{2,ac}(\R^n)=\W_2(\R^n) \cap \mathcal{P}_{ac}(\R^n)$.

There is a natural metric on $\W_2(\R^n)$, the quadratic Wasserstein metric, which is defined as
\begin{equation}\label{eq:W2-metric}
    W_2(\sigma,\mu) := \inf_{\pi\in\Gamma(\sigma,\mu)} \left(\int_{\R^{2n}}\|x-y\|^2d\pi(x,y)\right)^{\frac12},
\end{equation}
where $\Gamma(\sigma,\mu):=\{\gamma\in\mathcal{P}(\R^{2n}): \gamma(A\times \R^n) = \sigma(A),\; \gamma(\R^n\times A)=\mu(A) \textnormal{ for } A\subset\R^n\}$ denotes the set of couplings (measures on the product space with marginals $\sigma$ and $\mu$).

If $\sigma \in {\W_{2,ac}(\R^n)}$, then the following optimization problem has a solution
\begin{equation}\label{eq:W2_map}
    \min_{T:T_{\sharp}\sigma=\mu}\int_{\R^n}\|T(x)-x\|^2 \, d\sigma(x),
\end{equation}
where $T$ is a map in $L^2(\sigma)$ and $\sharp$ denotes the pushforward operation, $T_{\sharp}\sigma(A)=\sigma(T^{-1}(A))$ for $A$ measurable. Furthermore, the optimal coupling \eqref{eq:W2-metric} has the form $\pi = (\operatorname{id},T_{\sigma}^{\mu})_{\sharp}\sigma$, where $T_{\sigma}^{\mu}$ is the (up to additive constants, unique) solution to \eqref{eq:W2_map} \cite{brenier1991}. The Wasserstein metric can then be written as
\begin{equation}\label{eq:wass}
    W_2(\sigma,\mu)= \|T_{\sigma}^{\mu}-\operatorname{id}\|_{\sigma} :=\left( \int_{\R^n}\|T_{\sigma}^{\mu}(x)-x\|^2 \, d\sigma(x)\right)^{1/2}.
\end{equation}
 The optimal transport map $T_{\sigma}^{\mu}$ can be written as the gradient of a {convex function}, i.e.\ $T_{\sigma}^{\mu}=\nabla \varphi$ with $\varphi$ convex \cite{brenier1991}. We call a map $S$, which is the gradient of a convex function (but not necessarily the optimal transport map between two measures) a \emph{Brenier map}.

In the case of one-dimensional measures, the optimal transport map (and the Wasserstein distance) can be explicitly computed: For $\sigma,\mu \in \mathcal{P}(\R)$ and $\sigma \in\mathcal{P}_{ac}(\R)$, we get
\begin{equation}\label{eq:1-dOTmap}
    T_{\sigma}^{\mu} =  F_{\mu}^{-1}\circ F_{\sigma},
\end{equation}
where $F_{\sigma}: \R \to [0,1]$ is the cumulative distribution function (CDF) of $\sigma$, defined by
$
    F_{\sigma}(x)=\sigma((-\infty,x]).
$    
Here $F_{\mu}^{-1}$ denotes the pseudo-inverse $F_{\mu}^{-1}(y)=\min_z\{z:F_{\mu}(z)\geq y\}$.
Then we get
\begin{equation}\label{eq:W2-1d}
    W_2(\sigma,\mu) =\left( \int_{0}^1 |F_{\mu}^{-1}(x)-F_{\sigma}^{-1}(x)|^2\,dx\right)^{1/2}.
\end{equation}
Note that the assignment $\sigma \mapsto F_{\sigma}^{-1}$ is an isometry to the space $L^2(\R)$ with linear $L^2$ metric, {the geodesic of which is a straight line.}

In this paper, we deal both with probability measures on $\R^n$ and $\R$ and we denote the Wasserstein distance and optimal transport maps between such measures with the same symbol, independent of the dimension of the measures. It will be clear from context if the measures are on $\R^n$ or $\R$.

\subsection{Sliced Wasserstein distances}
The computation of the Wasserstein distance \eqref{eq:W2-metric} can be expensive, in particular in high dimensions (solving the linear program \eqref{eq:W2-metric} leads to $O(n^3\log(n))$). Therefore, there is interest in approximations of the Wasserstein distance which can be computed more efficiently. A well-studied class of approximations are the Sinkhorn distances \cite{cuturi-2013}, which add a regularization term to the linear program \eqref{eq:W2-metric}. The resulting problem can then be solved with matrix scaling algorithms \cite{sinkhorn67}.

In this paper we are interested in \emph{sliced Wasserstein distances}, which make use of the fact that the Wasserstein distance can be computed easily for one-dimensional measures, see \eqref{eq:W2-1d}. The main idea is to project $\sigma, \mu \in \mathcal{P}(\R^n)$ to a line parallel to  $\theta$, compute the one-dimensional distances between the projected measures, and then sum (or integrate) over all directions $\theta$. {More precisely, for}
 $\theta \in S^{n-1}$ {we define the projection}    $\mathcal{P}_{\theta}: \R^n \to \R$ by 
$\mathcal{P}_{\theta}(x) = x \cdot \theta = \left< x,\theta\right>$, {and denote} the projection of $\sigma \in \mathcal{P}(\R^n)$ by 
$\sigma^{\theta}:=(\mathcal{P}_{\theta})_{\sharp}\sigma$. Then the (continuous) sliced Wasserstein distance between $\sigma$ and $\mu$ is defined by
\begin{equation}\label{eq:cont-SW2}
    SW_2^2(\sigma,\mu) = \int_{S^{n-1}} W_2^2(\sigma^{\theta},\mu^{\theta}) \, du(\theta),
\end{equation}
where integration is over the uniform measure $u$ on $S^{n-1}$. Note that the Wasserstein distance under the integral is between one-dimensional measures, hence can be computed by \eqref{eq:W2-1d}. The sliced Wasserstein distance defines a metric {\cite{bonnotte13thesis}}. 

For a finite set $\{\theta_i\}_{i=1}^N$, we can consider a discrete version of \eqref{eq:cont-SW2}:
\begin{equation*}
   SW_2^2(\sigma,\mu) \approx \frac1N\sum_{i=1}^N W_2^2(\sigma^{\theta_i},\mu^{\theta_i}).
\end{equation*}
In practice, the sliced Wasserstein distance can be computed in this way by drawing i.i.d.\ samples $\theta_i$ from the uniform distribution on $S^{n-1}$.

We list some results on the relation between the sliced Wasserstein distance and the regular Wasserstein distance:
\begin{lemma}[\cite{bonnotte13thesis}]\label{lm:SWproperty}
Let $\theta \in S^{n-1}$ and $\sigma,\mu \in \mathcal{P}(\R^n)$ with $\sigma \in \mathcal{P}_{ac}(\R^n)$. Then we have
\begin{enumerate}
    \item \begin{equation*}
    W_2(\sigma^{\theta},\mu^{\theta})\leq \|\theta\cdot (id- T_{\sigma}^\mu)\|_{\sigma}\leq W_2(\sigma,\mu),
\end{equation*}
which implies $SW_2(\sigma,\mu)^2 \leq \frac{1}{n} W_2(\sigma,\mu)^2$.
\item If both $\sigma$ and $\mu$ have compact supports, then 
\begin{equation*}
    W_2(\sigma,\mu)^2\leq C_n SW_2(\sigma,\mu)^{\frac{1}{n+1}},
\end{equation*}
where $C_n$ is a positive constant depending on dimension $n$ and the supports of the measures.
\end{enumerate}
\end{lemma}
\begin{proof}

{See  Proposition 5.1.3 and Theorem 5.1.5 in \cite{bonnotte13thesis} for details.}
\end{proof}

Through the computation of Wasserstein distances between slices $\sigma^{\theta_i},\mu^{\theta_i}$, we also obtain one-dimensional optimal transport maps $\Tthetai$ as in \eqref{eq:1-dOTmap}. As suggested in \cite{pitie2007automated}, these can be stacked together to construct a map $\R^n \to \R^n$, which is not necessarily the optimal transport map between $\sigma$ and $\mu$, but can be used to define an iterative approximation scheme. This type of schemes, and in particular, their convergence, are the main topic of the present manuscript. We introduce them in \Cref{sec:slice-matching}.

\section{Slice-matching maps}\label{sec:slice-matching}
We are interested in the two types of slice-matching maps defined in \cite{pitie2007automated}. In this section we present a unifying framework which is closely related to \emph{compatible maps} \cite{khurana2022supervised,aldroubi20}.

\begin{definition}[Definition from \cite{pitie2007automated}]\label{def:single-matching}
Consider $\sigma \in \W_{2,ac}(\R^n)$, $\mu\in \W_2(\R^n)$.
\begin{enumerate}
    \item \textbf{Single-slice matching map:}  Let $\theta \in S^{n-1}$. The \newmap map from ${\sigma}$ to ${\mu}$ is defined by
\begin{equation}\label{eq:single-Sliceform}
 \singleslice(x) = x + (\Ttheta(x\cdot \theta)-x\cdot \theta)\,\theta,
\end{equation}
where $\Ttheta: \R\rightarrow\R$ denotes the optimal transport map from $\sigma^{\theta}$ to $\mu^{\theta}$.
\item \textbf{Matrix-slice matching map}: Let $\{\bth_1,\cdots, \bth_n\}\subset S^{n-1}$ be an orthonormal set. The \Mnewmap map from ${\sigma}$ to ${\mu}$ is defined by
\begin{equation}\label{eq:Sliceform}
   \Pslice(x) = x + P\begin{bmatrix}\Tthetaone (\bx\cdot \bth_1)-x\cdot \bth_1\\ \Tthetatwo(\bx\cdot \bth_2)-x\cdot \bth_2\\ \vdots\\ \Tthetan(\bx\cdot \bth_n)-\bx\cdot \bth_n
    \end{bmatrix},
\end{equation}
where $P = [\bth_1, \cdots, \bth_n]$ and $\Tthetai$ denotes the optimal transport map from $\sigma^{\theta_i}$ to $\mu^{\theta_i}$.
\end{enumerate}
\end{definition}
The following representation follows immediately from this definition:
$
          \Pslice(\bx) = \sum_{i=1}^n \Tthetai(\bx \cdot \bth_i)\, \bth_i.
$           
Similarly, if we choose an orthonormal set $\{\bth_1,\cdots, \bth_n\}\subset S^{n-1}$ with $\theta_1=\theta$, then
$
 \singleslice(\bx) =\Ttheta(\bx \cdot \bth)\theta + \sum_{i=2}^n (\bx \cdot \bth_i)\, \bth_i.
$    
This motivates the following unified framework:
\begin{definition}\label{def:unified-slicing}
    Let $\sigma \in \W_{2,ac}(\R^n)$, $\mu\in \W_2(\R^n)$. Choose $1\leq j \leq n$ and an orthonormal set $\{\bth_1,\cdots, \bth_n\}\subset S^{n-1}$. The $j$-slice matching map is defined by
    \begin{equation}\label{SliceOT-j}
         \Pjslice(\bx) = \sum_{i=1}^j \Tthetai(\bx \cdot \bth_i)\, \bth_i + \sum_{i=j+1}^n (\bx \cdot \bth_i)\, \bth_i,
\end{equation}
where $P = [\bth_1, \cdots, \bth_n]$.
\end{definition}

\begin{remark}\label{remark:single-slice}
We remark on some important properties of the $j$-slice matching map:
\begin{enumerate}
    \item Note that $T_{\sigma,P}^n =\Pslice$ and $T_{\sigma,P}^1 = {T_{\sigma,\theta_1}}$, where $\theta_1$ is the first column of $P$.
    \item  The map $\Pjslice$ is a Brenier map since it is the gradient of the convex function
    \begin{equation*}
        x\mapsto  \sum_{i=1}^j F_i(x\cdot \theta_i) + \sum_{i=j+1}^n \frac12 {(x\cdot \theta_i)^2},
    \end{equation*}
    where $F_i$ is the anti-derivative of $\Tthetai$. Therefore, $\Pjslice = T_{\sigma}^{\mu}$ if and only if $(\Pjslice)_{\sharp}\sigma = \mu$. 

    \item A computation shows that $\mathcal{P}_{\theta_{i}}\circ\Pjslice = \Tthetai \circ \mathcal{P}_{\theta_{i}}$ for $i=1,\ldots,j$ and $\mathcal{P}_{\theta_{i}}\circ\Pjslice = \mathcal{P}_{\theta_{i}}$ for $i=j+1,\ldots,n$.
    \item Note that (3) implies: If $\nu = (\Pjslice)_\sharp \sigma$, then $\nu^{\theta_i}=\mu^{\theta_i}$ for $1\leq i \leq j$. This property is crucial for the proof of our main result in \Cref{mainthm-Pj} and is the motivation for the name slice-matching map.
\end{enumerate}
\end{remark}

The map $\Pjslice$ is a special type of a \emph{compatible map}, as introduced in \cite{khurana2022supervised,aldroubi20}:
For a fixed $P \in O(n)$, the set of \textit{compatible maps} is defined by
 \begin{equation}\label{eq:compatibleSets}
     \mathfrak{S}(P)=\Big\{\bx \mapsto P \begin{bmatrix}
f_1\left(\left(P^t\bx\right)_1\right)\\
f_2\left(\left(P^t\bx\right)_2\right)\\
\vdots\\
f_n\left(\left(P^t\bx\right)_n\right)
\end{bmatrix}: f_i:\R\rightarrow\R ~\textrm{is increasing}  \Big\}.
 \end{equation}
Note that functions in this set can be written as $\sum_{i=1}^n f_i(\bx \cdot \bth_i)\, \bth_i$ with $P=[\theta_1,\ldots,\theta_n]$, which directly shows the relation to $\Pjslice$.

\begin{remark}
Direct verification shows that $ \mathfrak{S}(P)$ has a semi-group structure with composition as the group operation. {Note that, in general, the composition of two optimal transport maps is not necessarily an optimal transport map; however, in this case it is.} Moreover, if each $f_i$ {associated} with a compatible map given in \eqref{eq:compatibleSets} has an inverse, then the compatible map has an inverse in $ \mathfrak{S}(P)$, to which the corresponding $f_i$'s are simply replaced by their inverses. Note that the slice-matching maps of \Cref{def:unified-slicing} can be inverted easily.
\end{remark}

Compatible maps have been identified as a set of maps which allow for linear separability of two classes of measures in a tangent space, see \cite{khurana2022supervised,aldroubi20,park2018cumulative,moosmueller2020linear}. Results in these papers concerning supervised learning in the Wasserstein space therefore also hold for slice-matching maps.

\subsection{Iterative schemes via slice-matching maps}
Following \cite{pitie2007automated}, we define an iterative scheme using the slice-matching maps of \Cref{def:unified-slicing}:
\begin{definition}\label{def:stoch-iterative}
Let $\sigma_0\in \W_{2,ac}(\R^n)$ and $\mu\in \W_2(\R^n)$. For $k\geq 0$, choose $P_k \in O(n)$ and $\gamma_k \in [0,1]$. Fix $1\leq j \leq n$. Define
  \begin{equation}\label{sigmaiter-generalized}
      \sigma_{k+1}= \left((1-\gamma_k)\operatorname{id} +\gamma_k\Pkjslice\right)_{\sharp}\sigma_k, 
      \quad k\geq 0.
  \end{equation}
\end{definition}

\begin{remark}
Note that for $\gamma_k=1$ and $j=n$, we obtain the original scheme of \cite{pitie2007automated}, see \eqref{intro:iter}. We mention that our convergence results of \Cref{sec:iterative_scheme} do not hold for this scheme, since $\gamma_k=1$ does not satisfy \Cref{convergence-assumptions}. Results detailed in this section, however, hold for all schemes of \Cref{def:stoch-iterative}.
\end{remark}
We note that $(1-\gamma_k)\operatorname{id} +\gamma_k\Pkjslice$ is a Brenier map, since it is a convex combination of two gradients of convex functions. It is furthermore the optimal transport map from $\sigma_k$ to $\sigma_{k+1}$.

We show two results for the iterative scheme of \Cref{def:stoch-iterative}:
\begin{lemma}\label{lem:single-slice-stationary}
Let $\sigma_k \in \W_{2,ac}(\R^n)$ and $\mu \in \W_2(\R^n)$.  If $\gamma_k=1$ and $P_{k+1} = P_{k}$, then the iterative scheme \eqref{sigmaiter-generalized} gives $\sigma_{k+2} = \sigma_{k+1}$.
\end{lemma}
\begin{proof}
We first note that $\sigma_{k+1}^{\theta^k_{i}} = \mu^{\theta^k_{i}}$ where $P_k = [\theta_1^k,\ldots,\theta_n^k]$. Using \Cref{remark:single-slice}, we obtain
\begin{equation*}
    \sigma_{k+1}^{\theta^k_{i}} = (\mathcal{P}_{\theta^k_{i}} \circ \Pkjslice)_{\sharp}\sigma_{k} = \left(\Tkthetaik\right)_{\sharp}\sigma^{\theta^k_{i}}_k = \mu^{\theta^k_{i}}, \quad i=1,\ldots,j,
\end{equation*}
{where $\Tkthetaik$ denotes the optimal transport map between $\sigma^{\theta^k_{i}}_k$ and $\mu^{\theta^k_{i}}$.  }
This implies $\Tkplusthetaik = \operatorname{id}, i=1,\ldots,j,$ and therefore 
\begin{equation*}
  \Pkplusjslice(x) = \sum_{i=1}^j \Tkplusthetaik (x\cdot \theta^k_{i})\theta^k_{i} + \sum_{i=j+1}^n (x\cdot \theta^k_{i})\theta^k_{i} =x,
\end{equation*}
which implies $\sigma_{k+2} = \sigma_{k+1}$.
\end{proof}

\begin{remark}
    Note that for $P_{k+1}=P_k$ and $\gamma_k\neq 1$, the scheme does not necessarily become stationary:
    \begin{equation*}
        \sigma_{k+1}^{\theta^k_{i}} = \left((1-\gamma_k)\id + \gamma_k {\Tkthetaik}\right)_{\sharp}\sigma^{\theta^k_{i}}_k, \quad  i = 1,\ldots, j.
    \end{equation*}
    The iteration evaluates on the geodesic connecting the slices ${\sigma^{\theta^k_{i}}_k}$ and ${\mu^{\theta^k_{i}}}$, but is not necessarily equal to ${\mu^{\theta^k_{i}}}$.
\end{remark}

The following Lemma is crucial for the proof of our main result, \Cref{mainthm-Pj}, {as it relates} consecutive iterates of \eqref{sigmaiter-generalized} to a type of discrete sliced Wasserstein distance:

\begin{lemma}\label{lem:single-slice-distances}
Let $\sigma_k \in \W_{2,ac}(\R^n)$ and $\mu\in \W_2(\R^n)$. Let {$\sigma_{k+1}$} be defined through the iteration \eqref{sigmaiter-generalized}. Then
\begin{equation*}
    W_2^2(\sigma_{k+1},\sigma_{k}) = \gamma_k^2 \, \sum_{i=1}^jW_2^2(\sigma_{k}^{\theta_i^k},\mu^{\theta_i^k}),
\end{equation*}
where $P_k = [\theta_1^k,\ldots,\theta_n^k]\in O(n)$.
\end{lemma}

\begin{proof}

This result follows from direct computation:
\begin{align*}
    W_2(\sigma_{k+1},\sigma_{k})^2 &= \|T^{\sigma_{k+1}}_{\sigma_{k}} - \operatorname{id}\|_2^2 
    = \int_{\R^n} \|(1-\gamma_k)x +\gamma_k\Pkjslice(x) -x\|^2_2 \, d\sigma_k(x) \\
    & = \gamma_k^2 \, \int_{\R^n} \left\|\sum_{i=1}^j(\Tkthetaik(x\cdot \theta^k_i) -x\cdot \theta^k_i)\theta^k_i  \right\|_2^2 \, d\sigma_k(x) 
     = \gamma_k^2 \int_{\R^n} \sum_{i=1}^j(\Tkthetaik(x\cdot \theta^k_i) -x\cdot \theta^k_i)^2 \, d\sigma_k(x)\\
    & = \gamma_k^2 \sum_{i=1}^j \int_{\R} (\Tkthetaik(y) -y)^2 \, d\sigma_k^{\theta^k_i}(y)
     = 
    \gamma_k^2 \, \sum_{i=1}^jW_2(\sigma_{k}^{\theta_i^k},\mu^{\theta_i^k})^2,
\end{align*}
where $\Tkthetaik$ denotes the optimal transport map from $\sigma^{\theta_i^k}$ to $\mu^{\theta_i^k}$.
\end{proof}
\begin{remark}
Proposition 5.2.7 in \cite{bonnotte13thesis} provides an analogous result for the case of discrete measures with finite supports and $\gamma_k=1,j=n$.
\end{remark}

\section{Convergence of the stochastic iterative scheme}\label{sec:iterative_scheme}
To illustrate the idea of the iterative schemes \eqref{sigmaiter-generalized} as a stochastic gradient descent procedure in the 2-Wasserstein space of a certain functional, we first highlight the single-slice matching case. We then show a general interpretation for the $j$-slice matching scheme \eqref{sigmaiter-generalized}
as well as an a.s.\ convergence proof. To see the connection between the single-slice and matrix-slice schemes, we note the following relation between the uniform probability measure on the sphere $S^{n-1}$ and the Haar probability measure on the orthogonal group $O(n)$:

\begin{remark}\label{rmk:meas_relations}
    By an explicit geometric construction of the Haar measure (see e.g., \cite{meckes_2019}), a random orthogonal matrix $P$ distributed according to $u_n$ can be constructed by choosing the first column as $\theta_1$ a random vector according to the uniform probability measure $u$
on $S^{n-1}$, and then choosing the subsequent columns according to the surface area measures on subsets of $S^{n-1}$ that are orthogonal to the previous columns.  
\end{remark}

\subsection{Stochastic gradient descent interpretation}\label{sec:introIterschemes}

\subsubsection{Stochastic single-slice matching}\label{sec:SGD_single}
When $j=1$, i.e.,  a single directional vector is chosen at each iteration, by \Cref{rmk:meas_relations} {and \Cref{remark:single-slice}}, the iterative scheme in \Cref{def:stoch-iterative} becomes: Given $\theta_k\in S^{n-1}$ and $\gamma_k \in [0,1]$, 
  \begin{equation}\label{sigmaiter}
      \sigma_{k+1}= \left((1-\gamma_k)\operatorname{id} +\gamma_k\singleslicek\right)_{\sharp}\sigma_k, 
      \quad k\geq 0.
  \end{equation}

When each $\theta_k$ is chosen independently according to $u$, the uniform probability measure on $S^{n-1}$, the above iterative scheme \eqref{sigmaiter} can be interpreted {as} a stochastic gradient descent of the following functional minimization problem:
\begin{equation}\label{minLsigma}
    \min_{\sigma\in \W_{2,ac}(\R^n)}L(\sigma), 
\end{equation}
where 
\begin{equation}\label{eq:Lfunctional}
   {H}(\sigma):= \frac{1}{2}SW_2^2(\sigma,\mu) = \frac{1}{2} \int_{S^{n-1}}W_2^2(\sigma^{\bth}, \mu^{\bth})du(\bth).
\end{equation}
The unique minimizer is $\mu$ by the fact that $SW_2$ is a metric. Applying \Cref{lem:single-slice-distances} with $\gamma_k=1$ and $j=1$, one can equivalently write \eqref{eq:Lfunctional} as
\begin{equation}\label{eq:eqiv_Lfunctional}
   {H}(\sigma)= \frac{1}{2}\int_{S^{n-1}}W_2^2(\sigma,(\singleslice)_{\sharp}\sigma)du(\theta).
\end{equation}
The connection to the stochastic gradient decent method can be  observed by computing the formal  Fr\'echet derivative of $L$, following the ideas in  \cite{zemel2019frechet, ambrosio2008gradient}. 
We first note that for the functional $F(\sigma):= \frac{1}{2}W_2^2(\sigma,\mu)$ defined on $\W_2(\R^n)$ with a fixed $\mu\in \W_2(\R^n)$, the following differentiability property holds: For any $\sigma\in \W_{2,ac}(\R^n)$,
\begin{equation*}
    \lim_{\sigma_1\rightarrow\sigma}\frac{F(\sigma_1)-F(\sigma)-\left<\id-T_{\sigma}^\mu, T_{\sigma}^{\sigma^1}-\id\right>_{\sigma}}{W_2(\sigma_1,\sigma)}=0,
\end{equation*}
see \cite[Corollary 10.2.7]{ambrosio2008gradient}.
Then the Fr\'echet type derivative for $F$ at $\sigma$, denoted as $ F^{\prime}(\sigma)$ is a functional on the tangent space (see e.g. \cite[Definition 8.5.1]{ambrosio2008gradient} and \cite[p. 938]{zemel2019frechet})
\begin{equation*}
    \textrm{Tan}_{\sigma}:= \overline{\{\lambda(T-\id): T  \textrm{~ is a Brenier map~}; \lambda>0\}}^{L^2(\sigma)},
\end{equation*}
and is given by 

$
    F^{\prime}(\sigma)= \id - T_{\sigma}^{\mu},
$    

using the Riesz representation theorem on $L^2(\sigma)$. Relating $F$ and  \eqref{eq:eqiv_Lfunctional} and observing that $\singleslice$ is a Brenier map (see \Cref{remark:single-slice}), we define the formal Fr\'echet derivative \footnote{\textcolor{black}{\cite{backhoff2025stochastic} defines the formal Fr\'echet derivative corresponding to the barycenter problem in a similar spirit.}}
\begin{equation}\label{eq:LprimeSig}
   {H^{\prime}}(\sigma) (x) := \int_{S^{n-1}}\Big(x-\singleslice(x)\Big)du(\bth),
\end{equation}
where $T_{\sigma,\mu; \theta}$ is as defined in \eqref{eq:single-Sliceform}. Note that again the functional ${H^{\prime}}(\sigma)$ is identified as a function in $ L^2(\sigma)$ by the Riesz representation theorem. 
\begin{remark}
For the analysis in this manuscript, only a formal notion of Fr\'echet derivative as defined in \eqref{eq:LprimeSig} is needed. Under additional assumptions, e.g., if $\sigma,\mu\in \W_{2,ac}(K)$ where $K$ is a compact subset of $\R^n$, we get
\begin{equation*}
   \int_{K}{H^{\prime}}(\sigma)(x)\zeta(x)d\sigma(x) = \lim_{\epsilon\rightarrow0}
    \frac{SW_2^2(\left(\id+\epsilon \zeta\right)_{\sharp}\sigma,\mu)-W_2^2(\sigma,\mu)}{2\epsilon},
\end{equation*}
for any test diffeomorphism {$\zeta$} of $K$, see  \cite[Proposition 5.1.7]{bonnotte13thesis}.
\end{remark}

To develop an intuitive understanding of the scheme  \eqref{sigmaiter} as a stochastic gradient descent step for the minimization problem stated in \eqref{eq:Lfunctional}, we note that the stochastic version of Fr\'echet derivative \eqref{eq:LprimeSig} at a random  $\theta_k\in S^{n-1}$ is  $\id-\singleslicek$. Hence given a step size $\gamma_k\in [0,1]$, the corresponding push-forward map between measures is given by 
\begin{equation}
    x \mapsto x - \gamma_k (x-\singleslicek(x)),
\end{equation}
which gives  \eqref{sigmaiter}.

\subsubsection{Stochastic \texorpdfstring{$j$}{j}-slice matching}\label{subsec:stochastic-j-slice}
More generally, when the first $j$ columns of an orthogonal matrix are used in each iteration, and each $P_k$ is chosen independently according to $u_n$, the  Haar probability measure on the orthogonal group $O(n)$, the iterative $j$-slice matching scheme \eqref{def:stoch-iterative} can be viewed as a stochastic gradient descent scheme for the following functional minimization problem:

 \begin{equation}\label{eq:LOjfunctional}
    \min_{\sigma\in \W_{2,ac}(\R^n)}{H_j}(\sigma) = \frac{1}{2}\int_{O(n)}\sum_{i=1}^jW_2^2(\sigma^{\theta_i},\mu^{\theta_i})\,du_n(P),
\end{equation}
where $P = [\theta_1,\cdots,\theta_n]$. Applying \Cref{lem:single-slice-distances} with $\gamma_k=1$, one can equivalently write \Cref{eq:LOjfunctional} in the following way:
    \begin{equation}\label{eq:Lj-nWass}
       {H_j}(\sigma)=\frac{1}{2} \int_{O(n)}W_2^2(\sigma,(\Pjslice)_{\sharp}\sigma)\,du_n(P).
    \end{equation}
Following similar analysis in \Cref{sec:SGD_single}, we define the formal  Fr\'echet derivative
\begin{equation}\label{FrecheDer-Pj}
   {{H_j^{\prime}}}(\sigma)(x) = \int_{O(n)} (x-\Pjslice(x))\,du_n(P).
\end{equation}

\begin{remark}\label{rmk:Ljrelation}
\Cref{rmk:meas_relations} implies
\begin{enumerate}
    \item ${H_1 = H}$, as defined in \eqref{eq:Lfunctional},
    \item {${H_j}(\sigma)= \frac{j}{2}SW_2^2(\sigma,\mu)$} and hence ${H_j}(\sigma)=0$ if and only if $\sigma = \mu$.
    \item {$  {H_j^{\prime}}(\sigma)(x) = j \int_{S^n-1}\theta\left(\theta\cdot x - \Ttheta\left(\theta\cdot x\right)\right)du(\theta)$}, where $\Ttheta$ denotes the optimal transport map from $\sigma^{\theta}$ to $\mu^{\theta}$.

\end{enumerate}
\end{remark}

\subsection{Almost sure convergence of stochastic \texorpdfstring{$j$}{j}-slice matching}

For the convergence analysis of the iterative scheme defined in  \eqref{sigmaiter-generalized}, we need the following assumptions, which are adapted from \cite{backhoff2025stochastic}:

\begin{assumption}\label{convergence-assumptions}
Fix $\sigma_0,\mu \in \W_{2,ac}$ . Let $P_k\in O(n)$ and $\gamma_k \in [0,1]$ for $k\geq 0$.
\begin{enumerate}
   \item[(A1)] Given $\sigma_0,\mu \in {\W}_{2,ac}(\R^n)$, {there exists some compact set $K_{\sigma_0,\mu}\subseteq {\W}_{2,ac}(\R^n)$  such that:
  \begin{enumerate}
      \item[(i)] the sequence $\{\sigma_k\}_{k\geq 1}$ given in \eqref{sigmaiter} stays $K_{\sigma_0,\mu}$ independent of the choices of $\{\gamma_k\}_{k\geq 0}$ and $\{P_k\}_{k\geq 0}$. 
      \item[(ii)] ${H_j^{\prime}}(\eta) =0$ ($\eta$-a.e.) where   ${\eta}\in K_{\sigma_0,\mu}$ implies that $\eta = \mu$. \footnote{This condition is more generally referred to as the uniqueness of the Karcher mean, see \cite{backhoff2025stochastic, zemel2019frechet}.}
  \end{enumerate}  
     }
       
    \item[(A2)] Step-size assumption:
    \begin{align}
        \sum_{k=0}^{\infty}\gamma_k &= \infty,\label{gammaksum}\\
        \sum_{k=0}^{\infty}\gamma_k^2&<\infty.\label{gammaksqsum}
    \end{align}
\end{enumerate}
\end{assumption}
{The sliced-Wasserstein flow may not converge to the target measure in general \cite{cozzi2025long}. The authors prove convergence when the target measure is the Gaussian distribution and conjecture convergence under a finite entropy assumption on the source measure. As noted in Section 1.3, Assumption (A1-ii) is central and closely tied to the uniqueness of Karcher means in barycenter problems, the general understanding of which remains largely open. The study of sufficient conditions for (A1-ii) is ongoing work. In the remarks below, we present examples of sufficient conditions for (A1-i) as well as a narrow case for (A1-ii).}

\begin{remark}\label{compact_W2}
Observe that $ \{\sigma_k\}\subseteq {\W}_{2,ac}(\R^n)$,  given $\sigma_0,\mu \in {\W}_{2,ac}(\R^n)$, see \cite[Lemma C.10]{limoos2024approximation}. The existence of a compact set $K_{\sigma_0,\mu}\subseteq \W_{2}(\R^n)$ to which $\{\sigma_k\}$ belongs can be guaranteed 
with relatively mild conditions on $\sigma_0$ and $\mu$. One sufficient condition is that target measure $\mu$ has compact support, then $\{\sigma_k\}$ lies in ${\W_{2,ac}(\Omega)}$ for some compact set $\Omega\subseteq \R^n$, see \Cref{prop:A1_cond_1}. Note that $\W_{2}(\Omega)$ is compact with respect to the $W_2$ metric by tightness and \cite[Theorem 5.9]{santambrogio2015optimal}. When $\mu$ does not have compact support, for example when $\mu$ is the Gaussian,  a compact subset can be constructed in $W_2(\R^n)$ given that $\sigma_0,\mu$ has finite $p$-th moments for $p>2$.  For example, if $\sigma_0$ and $\mu$ have finite third moments, then $M_3(\sigma_k)\leq M$ for some positive constant $M$ depending on the third moments of $\sigma_0$ and $\mu$, see {\Cref{prop:uniform3rdMoment}}. Note that the closure of $\{\nu\in \W_{2,ac}(\R^n): \int_{\R^n}\|x\|^3d\nu(x)\leq M\}$ is compact in $\W_{2}(\R^n)$ \cite{cozzi2025long}, see also {\Cref{prop:MomentSetCpt}}.
\end{remark}

\begin{remark}
 
 To {fulfill} (A1-i),  $K_{\sigma_0,\mu}$ needs to be a compact subset of absolutely continuous measures in $\W_{2}(\R^n)$, i.e.,  $K_{\sigma_0,\mu}\subseteq \W_{2,ac}(\R^n)$. One sufficient condition on which any limit point of $\{\sigma_k\}_{k\geq 1}$ stays absolutely continuous is that the densities of probability measures 
 $\sigma_k$ satisfy the energy bounds in \eqref{finite_G_energy}, where $G$ is a lower semi-continuous function and convex function on $[0,\infty)$ such that \eqref{G_cond} holds, see \Cref{lm:abs_continuity_energy}. For example, for $G(t)=t\log(t)$, $K_{\sigma_0,\mu}$ can be constructed as in \Cref{compact_W2} with an additional constraint that  the measures have uniformly bounded entropy.

\end{remark}

{\begin{remark} A sufficient scenario for {(A1-ii) and $j=1$} can be found in \cite[Lemma 5.7.2]{bonnotte13thesis}; specifically,  if the $K_{\sigma_0,\mu}$  {satisfying (A1-i}) is a  subset of $\W_{2,ac}(B_r)$ with strictly positive density, {then $H_1^{\prime}(\eta) = 0$ ($\eta$-a.e.)  for $\eta \in  K_{\sigma_0,\mu}$  implies that $\eta = \mu$}. Here $B_r\subseteq \R^n$ is the  closed ball of radius $r$ centered at the origin. { Note that the same condition works for all  $1\leq j\leq n$ since $H_j^\prime = j H_1^\prime$ (see part (3) of \Cref{rmk:Ljrelation}). }
\end{remark}

\begin{remark}
    Choices of $\gamma_k$ satisfying (A2) include $\gamma_k = \frac{1}{k}, k\geq 1$ and $\gamma_k = \frac{1+\log_2(k)}{k}$. Note that $\gamma_k=1$ does not satisfy the assumptions and therefore, our convergence result of \Cref{mainthm-Pj} does not hold for the original scheme of \cite{pitie2007automated}.
\end{remark}

\begin{theorem}[Main result]\label{mainthm-Pj}
    Let $\sigma_0,\mu\in \W_{2,ac}(\R^n)$ and $P_k \overset{\textrm{i.i.d}}{\sim} u_n, k\geq 0$,  where $u_n$ is the Haar probability measure on $O(n)$. Fix $1\leq j \leq n$ and assume that {(A1) and (A2)} hold for ${H_j}$ of \eqref{eq:LOjfunctional}.  Then the $j$-slice matching scheme \eqref{def:stoch-iterative} satisfies
    \begin{equation*}
         \sigma_k \xrightarrow{W_2} \mu, \quad \text{a.s.\ as } k \to \infty.
    \end{equation*} 
\end{theorem}

\begin{proof}\label{proof:mainthm}
        The proof is based on a careful modification of the proof of Theorem 1.4 in \cite{backhoff2025stochastic}, a result by Backhoff-Veragyas et al.\ solving a  population  barycenter problem through
\begin{equation*}
    \min_{\sigma\in \W_{2,ac}(\R^n) }\mathcal{F}(\sigma)= \frac{1}{2}\int_{\W_{2,ac}(\R^n)} W_2^2(\sigma,m)\,d\Pi(m),
\end{equation*}
    where $\Pi$ is a probability measure defined on the space $\W_2(\R^n)$ of probability measures, which gives full measure to a compact subset.
In the present manuscript, we seek to recover $\mu$ through 
\begin{equation*}
   \min\limits_{\sigma\in \W_{2,ac}(\R^n)}{H_j}(\sigma)=\frac{1}{2}\int_{O(n)}W_2^2(\sigma, (\Pjslice)_{\sharp}\sigma)\,du_n(P), 
\end{equation*}
see \Cref{subsec:stochastic-j-slice}. Note that by the change-of-variable formula, ${H_j}$ becomes very similar to $\mathcal{F}$:
\begin{equation*}
   {H_j}(\sigma)= \frac12\int_{\varphi_{\sigma,\mu}^j(O(n))}W_2^2(\sigma,m) \, d\Pi_{\sigma,\mu}^j(m)
\end{equation*}
where $\varphi_{\sigma,\mu}^j:O(n) \to \mathcal{W}_2(\R^n)$ is defined by $\varphi_{\sigma,\mu}^j(P) = (\Pjslice)_{\sharp}\sigma$ and $\Pi_{\sigma,\mu}^j = {\varphi_{\sigma,\mu}^j}_{\sharp}u_n$. Through this interpretation, we arrive at a barycenter problem on a subset of $\mathcal{W}_2(\R^n)$ which is parametrized by $O(n)$.

Note that $\varphi_{\sigma,\mu}^j(O(n))$ is a compact subset to which $\Pi_{\sigma,\mu}^j$ gives full measure, however, the key difference to $\mathcal{F}$ lies in the dependence of both $\varphi_{\sigma,\mu}^j(O(n))$ and $\Pi_{\sigma,\mu}^j$ on $\sigma$. Therefore, the proof of \cite{backhoff2025stochastic} does not directly carry over, and it is more natural to work with integration over $O(n)$ than over $\Pi_{\sigma,\mu}^j$.

{We follow the outline of the proof of  \cite[Theorem 1]{backhoff2025stochastic} and show a careful adaptation. We reproduce relevant details for completeness.} 
The key modification is to show analogs of Proposition 3.2 and Lemma 3.3 in \cite{backhoff2025stochastic}, which are given by \Cref{analog3.2,analog3.3} in \Cref{appendixProofs}. In essence, \Cref{analog3.2} shows that the sequence $\{{H_j}(\sigma_k)\}_{k\in \mathbb{N}}$ is decreasing in expectation, and \Cref{analog3.3} verifies continuity properties of  ${H_j}$ and ${H_j^{\prime}}$.
{More specifically, we first
note that the unique minimizer for \eqref{minLsigma} is $\mu$ and ${H_j}(\mu)=0$, and  introduce similarly $l_i: ={H_j}(\sigma_i)$ and $\alpha_i:= \Pi_{k=1}^{i-1}1/(1+\gamma_k^2)$. By \eqref{gammaksqsum} in (A2), the sequence $(\alpha_i)$ converges to some finite $\alpha_{\infty}>0$. Then by \Cref{analog3.2}
\begin{equation*}
  \mathbb{E}[l_{i+1}-(1+\gamma_i^2)l_i| \mathcal{F}_i] \leq \gamma_i^2 H_j(\mu)- \gamma_i\|H_j^\prime (\sigma_i)\|^2_{\sigma_i}\leq 0,
\end{equation*}
where $ \mathcal{F}_i$ is the associated $\sigma$-algebra generated by $P_0,...,P_i$, see \Cref{Pkdef}. 
Multiplying by $\alpha_{i+1}$, one similarly obtains
\begin{equation}\label{martingale}
    \mathbb{E}[{\alpha_{i+1}l_{i+1} - \alpha_i l_i| \mathcal{F}_i}]\leq \alpha_{i+1}\gamma_i^2H_j(\mu)-\alpha_{i+1}\|H_j^\prime(\sigma_i)\|^2_{\sigma_i}\leq 0,
\end{equation}
which implies $\mathbb{E}[\widehat{l}_{i+1}-\widehat{l}_i|\mathcal{F}_i]\leq 0$, where $\widehat{l}_i := \alpha_i l_i$. In other words, $(\widehat{l}_i)_{i\geq 0}$ is a supermartingale with respect to $(\mathcal{F}_i)$. It is not hard to see that $\widehat{l}_i$ is uniformly bounded from below, and hence by the supermartingale convergence theorem \cite[Corollary 11.7]{williams1991probability} we have that $\widehat{l}_i\rightarrow \widehat{l}_{\infty}$ a.s. (and hence $l_i \rightarrow l_{\infty}$) for some non-negative random variable $l_{\infty}\in L^1$. To see that $l_\infty=0$, we take the expectations in \eqref{martingale} and sum over $i$ using telescoping to obtain 
$\mathbb{E}[\alpha_{i+1}l_{i+1}] - \mathbb{E}[l_0\alpha_0]\leq -\sum\limits_{k=1}^i \alpha_{k+1}\gamma_k\mathbb{E}[\|H_j^\prime (\sigma_k)\|^2_{\sigma_k}]$. Then, taking the liminf, we have
$-\infty < \mathbb{E}[\alpha_{\infty}l_{\infty}]-\mathbb{E}[\alpha_{0}l_{0}] \leq - \mathbb{E}[\sum\limits_{k=1}^{\infty} \alpha_{k+1}\gamma_k\|H_j^\prime (\sigma_k)\|^2_{\sigma_k}]$, where the left inequality follows from Fatou and the right inequality follows from monotone convergence. Since $\lim_{k\rightarrow \infty} \alpha_k >0$, it follows that $\sum\limits_{i=1}^{\infty} \gamma_i\|H_j^\prime (\sigma_i)\|^2_{\sigma_i}<\infty \quad $ a.s., which, combined with \eqref{gammaksum} in (A2), further implies that 
\begin{equation}\label{eq:zeroderlimit}
\liminf_{i\rightarrow \infty} \|{H_j^\prime}(\sigma_i)\|^2_{\sigma_i}=0 \quad \textrm{a.s.}
\end{equation}
Now {let $K_{\sigma_0,\mu}$ be a compact set such that $\sigma_i\in K_{\sigma_0,\mu}$ for all $i$} (see (A1-i)).
Let $\epsilon>0$ and $K_{\epsilon}= \{\sigma :{H_j}(\sigma)\geq\epsilon\}\cap K_{\sigma,\mu}$. 
Since $K_\epsilon$ is compact, we have $\inf \limits_{K_\epsilon}\|{H_j}^{\prime}(\sigma)\|_{\sigma}>0$. Otherwise, it would imply the existence of $\sigma_{\epsilon} \in K_{\epsilon}$ such that ${H_j^{\prime}}(\sigma_{\epsilon}) =0$ by part (ii) of Lemma \ref{analog3.3}. However, this would lead to contradictory statements: $\sigma_{\epsilon}=\mu$ according to (A2-ii), and ${H_j}(\sigma_\epsilon)\geq \epsilon$ by part (i) of Lemma \ref{analog3.3}. Using $\inf \limits_{K_\epsilon}\|{H_j}^{\prime}(\sigma)\|_{\sigma}>0$, we deduce similar relationships between events as in \cite{backhoff2025stochastic}: 
$ \{l_\infty \geq 2\epsilon\} \subset \left\{ \sigma_i \in K_{\epsilon} \textrm{~for all $i$ large enough}\right\} \subset \bigcup_{m\in \mathbb{N}} \left\{ \|H_j^\prime(\sigma_i)\|^2_{\sigma_i}> 1/m: \textrm{~ for all $i$ large enough}\right\}
$, and hence $\{l_\infty \geq 2\epsilon\} \subset \left\{ \liminf\limits_{i\rightarrow \infty} \|H_j^\prime(\sigma_i)\|^2_{\sigma_i}>0 \right\} $.
By \eqref{eq:zeroderlimit}, we have that 
  $l_{\infty} =0$ a.s. Then by the compactness of $K_{\sigma,\mu}$ and \Cref{analog3.3}, we have that $\sigma_i \xrightarrow{W_2} \mu$ a.s., which trivially implies that $SW_2(\sigma_i,\mu)\rightarrow 0$ a.s.\ since $SW_2\leq W_2$ (see e.g., \Cref{lm:SWproperty}). Hence we also have that $\sigma_i\xrightarrow{SW_2} \mu$}.
\end{proof}

{As a corollary of the relation between the sliced-Wasserstein and Wasserstein distances (see \Cref{lm:SWproperty}), the above $ \sigma_k \xrightarrow{SW_2} \mu$. } We also note that \Cref{mainthm-Pj} with $j=n$ shows a.s.\ convergence of the matrix-slice matching scheme \eqref{intro:randomized_iter}. For $j=1$ we also get convergence for the single-slice scheme, which is summarized in the following:
\begin{corollary}\label{mainthm}
      Let $\sigma_0,\mu\in \W_{2,ac}(\R^n)$ and $\bth_k \overset{\textrm{i.i.d}}{\sim} u, k\geq 0$,  where $u$ is the uniform probability measure on $S^{n-1}$.  Assume that (A1), (A2), (A3) hold for $L$ of \eqref{eq:Lfunctional}. Then the single-slice matching scheme \eqref{sigmaiter} satisfies
    \begin{equation*}
       \sigma_k \xrightarrow{W_2} \mu, \quad \text{a.s.\ as } k \to \infty.
    \end{equation*}   
\end{corollary}
\begin{proof}
The result follows directly from combining \Cref{mainthm-Pj} with \Cref{rmk:Ljrelation}.  
\end{proof}

We close this section with an example on translations.
The above theorem offers a straightforward application to verify the almost sure convergence of the stochastic single-slice scheme \eqref{sigmaiter} when the initial and target measure are related through a translation.
\begin{example}\label{ex:translation}
 Let $\sigma \in \W_{2,ac}(\R^n)$ and let $\mu = {T_b}_{\sharp}\sigma$ with $T_b(x)= x+b, b\in \R^n, b\neq 0$. We apply the single-slice iterative scheme (\eqref{sigmaiter-generalized} with $j=1$ and $\bth_k \overset{\textrm{i.i.d}}{\sim} u, k\geq 0$ on $S^{n-1}$) with $\sigma_0=\sigma$. Then for any choice of sequence $\gamma_k$ satisfying (A2), we get $\sigma_k\rightarrow \mu$ a.s. in $W_2$.   
\end{example}
\begin{proof}
    By \Cref{mainthm}, it suffices to show that the assumptions (A1) and  (A3) are satisfied.

\begin{enumerate}
   \item[{(A1-i)}] Let $K_{\mu}^b= \{\nu_z = (T_z)_{\sharp}\mu: T_z(x)=x+z, \|z\|\leq \|b\| \}$. The compactness of $ K_{\mu}^b$ can be seen via  the compactness of the ball $\{z\in\R^n: \|z\|\leq \|b\| \}$ and the continuity of the operator from $(\R^n,\|\cdot\|)$ to $(\mathcal{P}(\R^n),W_2)$, given by

    $
        z\mapsto (T_z)_{\sharp}\mu.
    $        

    Note that $\sigma \in K_{\mu}^b$.
    The fact that $\{\sigma_k\}$ stays in $ K_{\mu}^b$ can be seen from the geodesic convexity of $ K_{\mu}^b$, i.e., for any $\gamma\in [0,1]$ and $\nu_z\in  K_{\mu}^b$ we have

    $
        \left((1-\gamma)\id + \gamma {T_{\nu_z,\theta}}\right)_{\sharp}\nu_z \in  K_{\mu}^b.
    $        

    This follows from ${T_{\nu_z,\theta}}(x) = x-(\theta\cdot z)\theta$ and from 
 
    $
    (1-\gamma)x + \gamma {T_{\nu_z,\theta}}(x) = x - \gamma (\theta\cdot z)\theta.    
    $

    Therefore
    $\|\gamma(\theta\cdot z)\theta\|\leq \|z\| \leq \|b\|$ and $\left((1-\gamma)\id + \gamma {T_{\nu_z,\theta}}\right)_{\sharp}\nu_z \in  K_{\mu}^b$.

\item[{(A1-ii)}]
We will show that for any $\nu_z\in K_{\mu}^b$, 
$
 {H^{\prime}}(\nu_z)= 0 $ if and only if  $\nu_z= \mu.$  
Let $\nu_z = (T_z)_{\sharp}\mu$, where $T_z(x)=x+z$. A direct computation gives 
\begin{equation*}
   {H^{\prime}}(\nu_z)(x)= \int_{S^{n-1}} (\theta\cdot z)\theta\,du(\theta).
\end{equation*}
Let $Az := \int_{S^{n-1}} (\theta\cdot z)\theta\, du(\theta)$ where $A\in \R^{n\times n}$ and observe that for the standard basis vector $e_i\in \R^n$,
\begin{equation*}
   Ae_i= \int_{S^{n-1}} (\theta\cdot e_i)\theta\,du(\theta) = [0,\cdots,0, w_i, 0\cdots, 0]^t,
\end{equation*}
where all but the  $i$-th entry of the RHS vector are zero and  $w_i>0$ (see \Cref{sphereint_example}).  Since $A$ is a diagonal matrix with positive diagonal entries, it follows that  $Az = 0$ if and only if $z= 0$. \qedhere
\end{enumerate}
    \end{proof}

\begin{remark}\label{rem:convergence-rate}
For the $\sigma,\mu$ in \Cref{ex:translation}, the iterated scheme also  converges a.s. when $\gamma_k=1$. Moreover by \cite[Proposition 28]{limoos2024approximation} and a direct computation, we have $\mathbb{E}W_2^2(\sigma_k,\mu)= (\frac{n-j}{n})^k\|b\|^2$. In particular, when $j=n$, the target $\mu$ is matched in a single step, i.e., $\sigma_1 = \mu$. 
 
\end{remark}


\section{Numerical experiments}\label{sec:experiments}

{
In this section, we provide a series of numerical experiments illustrating the convergence properties of the slice-matching scheme. We consider imaging examples (continuous case) and Gaussians (discrete case) to study the influence of resolution, convergence rates, and dimensionality.

\subsection{Convergence of the Slice-Matching Scheme}\label{sec:matrix-num-convergence}
We first demonstrate the convergence behavior of the proposed scheme in an imaging context. As a representative example, we consider digit morphing using MNIST data. Specifically, we initialize with an image of the digit ``5'' and gradually morph it into the digit ``1'' by iteratively applying the slice-matching procedure. This experiment provides a visual illustration of the convergence dynamics, highlighting how successive iterations progressively align structural features of the two images. 

We show a numerical experiment morphing a digit $5$ into a digit $1$ using  the iterative scheme \eqref{sigmaiter-generalized} with $j=n=2$ and choice $\gamma_k = \frac{1+\log_2(k)}{k}$. The random orthogonal matrix 
\begin{equation*}
P = 
\begin{bmatrix}
  \cos\beta & \sin\beta\\-\sin\beta &  \cos\beta  
\end{bmatrix}    
\end{equation*}
is generated via choosing a uniformly distributed random angle $\beta \in [0,\pi/2)$. 

A subset of the first $25$ iterations are shown in \Cref{fig:digit_morph_10frames}. \Cref{fig:digit_morph_matrix_log} shows how  the relative sliced Wasserstein error $\frac{SW_2(\sigma_k,\mu)}{SW_2(\sigma_0,\mu)}$ decreases as the iteration variable $k$ increases. The first step $\sigma_1$ creates the biggest drop, as it translates the digit $5$ into the correct location. The other iteration steps ``stretch'' the $5$ into the $1$. There is some fluctuation in the convergence $SW_2(\sigma_k,\mu) \to 0$ since this is a stochastic iteration. There are still some artifacts present in $\sigma_{25}$ which are due to a combination of small number of iterations and numerical errors. The average of $\frac{SW_2(\sigma_9,\mu)}{SW_2(\sigma_0,\mu)}$  across 10 trials is $\approx 9.05 \times 10^{-3}$.

The fact that the error does not drop to a value closer to $0$, as expected from our theoretical results, has to do with numerical errors from the image resolution, as we discuss in \Cref{sec:image-resolution}.

We provide a similar numerical study using the single-slice scheme (rather than the matrix scheme) in \Cref{sec:additional-numerics}.


\begin{figure*}[t!]
\centering
\includegraphics[width=\linewidth]{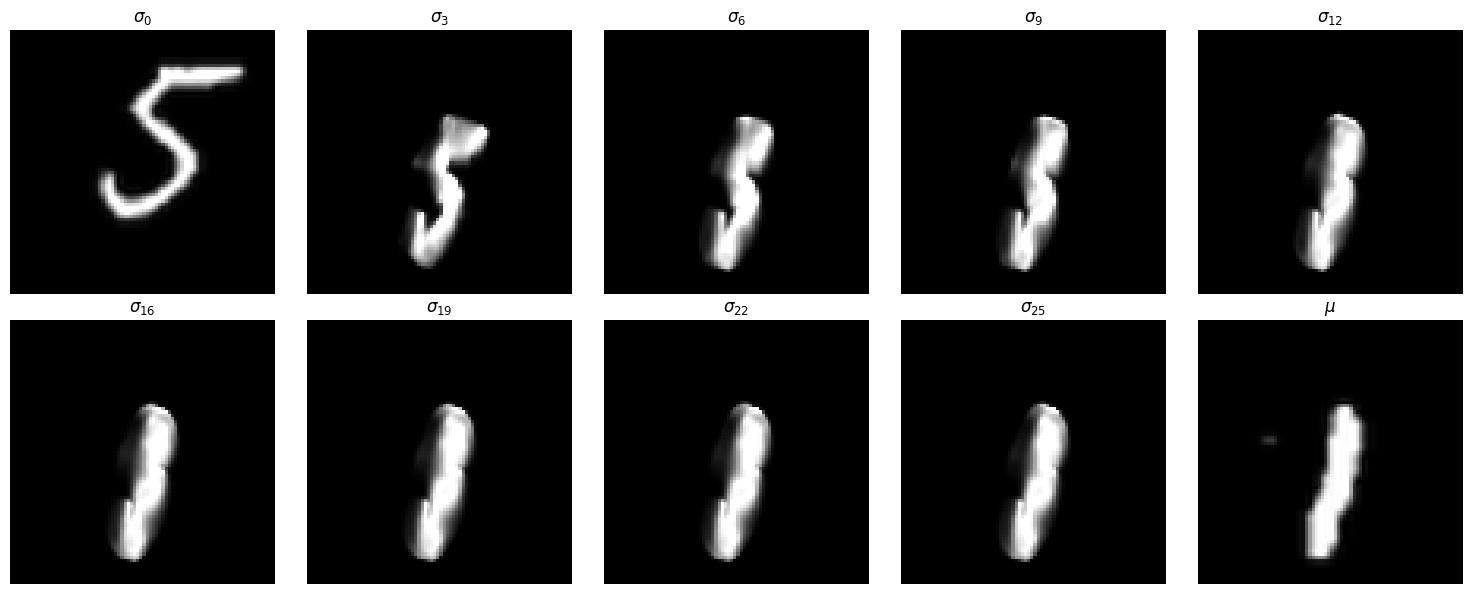}
\caption{Digit morphing experiment ($5 \to 1$), matrix slice version with $\gamma_k=\frac{1+\log_2(k+2)}{k+2}$. 
Ten equally spaced iterations $\sigma_0,\sigma_3,\sigma_6,\ldots,\mu$ are shown, capturing the progressive transformation.}
\label{fig:digit_morph_10frames}
\end{figure*}

\begin{figure*}[t!]
\centering
\includegraphics[scale=0.5]{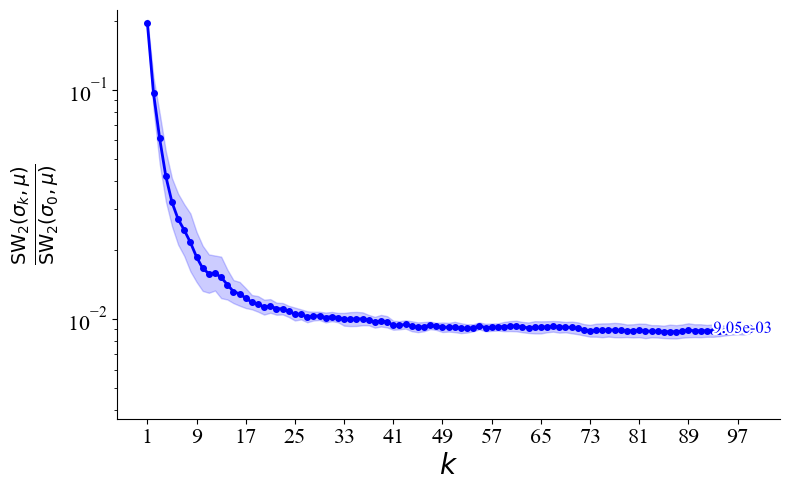}
\caption{Digit morphing experiment ($5 \to 1$) at $84 \times 84$ image resolution, matrix slice version with 
log-decaying step size $\gamma_k = \frac{1+\log_2(k+2)}{k+2}$. 
Intermediate morphs illustrate how the distribution gradually transforms from the source digit “5” to the target digit “1”. The figure shows the mean and standard deviation across 10 trials.}
\label{fig:digit_morph_matrix_log}
\end{figure*}

\subsection{Influence of Resolution on Convergence} \label{sec:image-resolution}
In the above experiment (\Cref{fig:digit_morph_10frames} and \Cref{fig:digit_morph_matrix_log}), we observe that the error curve does not converge to zero but rather stabilizes at a positive value. To further investigate this phenomenon, we conducted a systematic study on the influence of resolution and sampling. For the MNIST data (digits 5 and 1), we considered images at increasing resolutions ($(42 \times 42), (84 \times 84), \text{ and } (168 \times 168))$, see \Cref{fig:digits_resolution_log}. The figure shows how the error decreases with increasing image resolution.

\begin{figure*}[t!]
\centering
\includegraphics[width=0.9\linewidth]{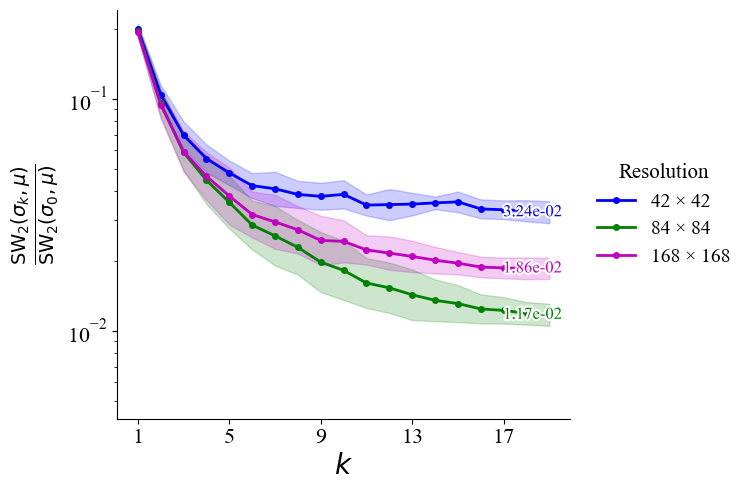}
\caption{Digit morphing experiment ($5 \to 1$) with sliced Wasserstein error (SW$_2$) under $\gamma_k=\frac{1+\log_2(k+2)}{k+2}$, 
for image resolutions $(42\times 42)$, $(84\times 84)$, and $(168\times 168)$. The curves show mean values and standard deviations over 10 trials.
The log-scaled y-axis highlights how finer image resolutions reduce the error and improve convergence.}
\label{fig:digits_resolution_log}
\end{figure*}

In addition, we merge an image of a 2D Gaussian distribution into its shifted version with varying resolution (\Cref{{fig:gaussian_shifted_SW_const}}) and the matrix-slicing scheme $(j=n=2)$. In this case, convergence is theoretically guaranteed after 1 step when choosing $(\gamma = 1)$, see \cite[Proposition 28]{limoos2024approximation}. This is why the curves are essentially flat as $k$ increases.
The numerical experiment shows that convergence improves with increased resolution, however, numerical errors prevent the scheme from dropping to a value much closer to $0$.

\begin{figure*}[t!]
\centering
\includegraphics[width=0.9\linewidth]{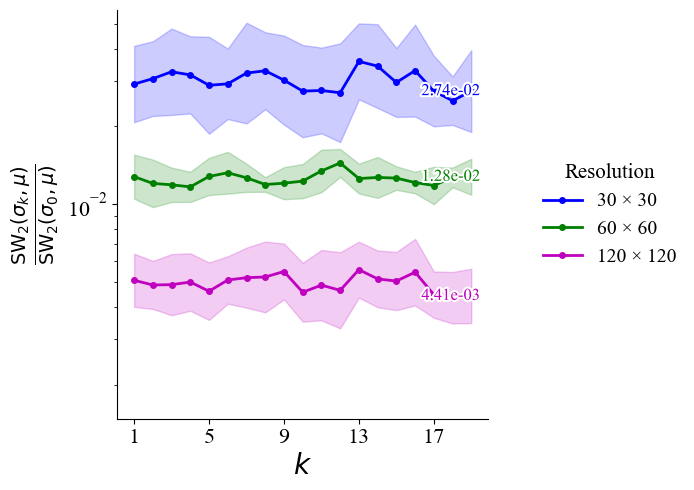}
\caption{Merging an image of a 2D Gaussian into its shifted version. The $y-$axis describes the sliced Wasserstein error (SW$_2$) for constant step size $\gamma=1$, and the $x-$axis describes the number of iterations.
Gaussian images are created with varying resolutions ($(30\times 30)$, $(60\times 60)$, and $(120\times 120)$). The curves show the mean value with standard deviation across 10 trials.
The plots shows how the error decreases with increased resolution.
}
\label{fig:gaussian_shifted_SW_const}
\end{figure*}

The results reveal that higher image resolution leads to improved convergence performance. For this type of image morphing task, the results suggest that ``perfect'' resolution would recover the expected convergence behavior. Thus, these findings emphasize the importance of discretization effects in practical implementations.

\subsection{Numerical Study of Convergence Rates}\label{sec:convergence-rate}
Our theoretical results establish convergence of the slice-matching scheme under certain assumptions, see \Cref{mainthm-Pj}. In specific cases, such as matching shifted measures, we are able to characterize the convergence rate analytically as mentioned in \Cref{rem:convergence-rate}. In particular, for a measure that is matched to its shifted version (shift by $b \in \mathbb{R}^n$ using $j$ slices and $\gamma_k=1$, the convergence goes as $\mathbb{E}W_2^2(\sigma_k,\mu)= (\frac{n-j}{n})^k\|b\|^2$).

\begin{figure*}[h!]
\centering
\includegraphics[width=0.9\linewidth]{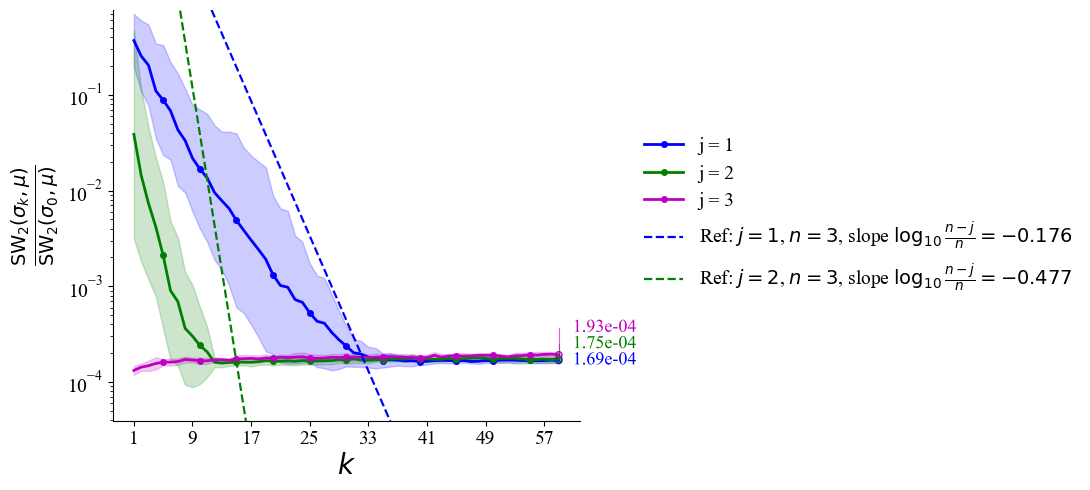}
\caption{Convergence rate for merging a Gaussian into its shifted version with $n=3$ and $j=1,2,3$, under $\gamma=1$. Number of discrete points are $N=2000$.
The plots show the relative sliced Wasserstein error $\frac{\mathrm{SW}_2(\sigma_k,\mu)}{\mathrm{SW}_2(\sigma_0,\mu)}$ on a log scale across iterations $k$. 
Shaded regions represent one standard deviation over 10 trials. 
Theoretical reference lines with slopes are included for $j=1,2$ to compare with the empirical decay rate. 
Note that for $j=3=n$, the error curve is flat, because the scheme converges after 1 step.  The error should theoretically vanish, but with sample size ($N=2000$), numerical fluctuations prevent it from reaching exact zero. }
\label{fig:conv_n3}
\end{figure*}

To validate this convergence rate numerically, we applied the slice-matching scheme to shifted Gaussian distributions and studied the decay of the error across iterations, see \Cref{fig:conv_n3}. Here we consider discrete Gaussians in 3D sampled with $N=2000$ rather than images of Gaussians. We compare the slice matching scheme using different amount of slices $j=1,2,3$, as the convergence rate depends on $j$. The figure shows both the relative error in sliced Wasserstein distance as the iteration variable $k$ increases (mean and standard deviation across 10 trials), and references lines (dashed) that indicate the theoretically established convergence rate. We see that the convergence behavior is very similar to the theoretical rate. Also, note that for $j=n=3$, the scheme converges after 1 step, this is why the purple line is flat.


We then examined whether the observed convergence rate extends to $\gamma_k = \frac{1+\log_2(k+2)}{k+2}$, for which we have convergence guarantees, but no theoretical proof for the convergence rates. Preliminary results suggest that the convergence rate is slower than for $\gamma=1$, see \Cref{fig:conv_n3_gammaComp}, although further theoretical investigation is required to fully establish generality.

\begin{figure*}[t!]
\centering
\includegraphics[width=0.9\linewidth]{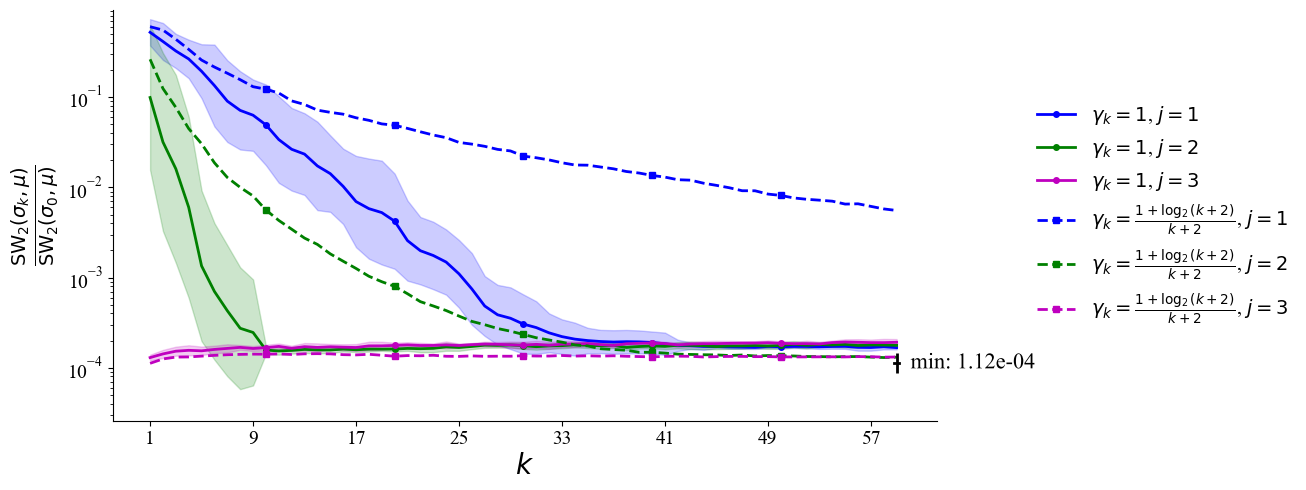}
\caption{Convergence rate for merging a Gaussian into its shifted version with $n=3$ and $j=1,2,3$, comparing constant step size ($\gamma=1$, solid lines) and log-decaying step size 
($\gamma_k = \frac{1+\log_2(k+2)}{k+2}$, dashed line). The number of sampling points for the Gaussian is $N=2000$.
The curves show the sliced Wasserstein error decays on a log scale, averaged over 10 trials with shaded error bands (for the $\gamma=1$ case; for the log-decaying step size, we only show the mean for better visualization). 
The comparison shows that the adaptive $\gamma_k$ scheme decays more slowly than the constant step size case, especially for smaller $j$. Again, we obtain a flat curve for $j=n=3$ because the scheme converges within 1 step.  The error should theoretically vanish, but with sample size ($N=2000$), numerical fluctuations prevent it from reaching exact zero. }
\label{fig:conv_n3_gammaComp}
\end{figure*}

\subsection{Higher-Dimensional Example}\label{sec:high-dim-examples}
Finally, we explored the influence of dimensionality on convergence. For this purpose, we considered matching shifted Gaussian distributions in higher dimensions, specifically in dimensions $(n = 15)$ and $(n = 100)$. In these cases, the slice-matching scheme again converges as expected (with appropriate rate as in \Cref{sec:convergence-rate} for $\gamma=1$), thereby providing empirical evidence that the method remains effective in moderately high-dimensional settings. These results suggest promising scalability of the approach, see \Cref{fig:conv_n100} for the dimension $n=100$.

\begin{figure*}[t!]
\centering
\includegraphics[width=0.9\linewidth]{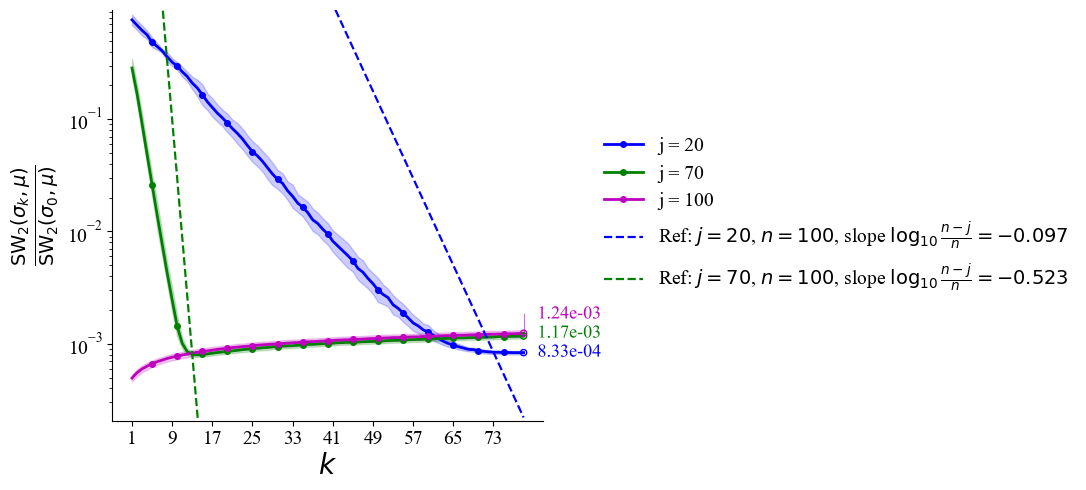}
\caption{The convergence rate for merging a Gaussian to its shifted version in dimension $n=100$ and  $j=20,70,100$, under step-size $\gamma=1$. 
The curves demonstrate how larger $j$ improve the convergence rate. The number of discrete points sampling points for the Gaussian is $N=2000$.
For $j=n=100$, the curve is flat because the scheme converges after 1 step. The error should theoretically vanish, but with sample size ($N=2000$), numerical fluctuations prevent it from reaching exact zero. }
\label{fig:conv_n100}
\end{figure*}

}

\section{Conclusion}
Motivated by the availability of closed-form formulae for one-dimensional optimal transport maps and the associated computational advantages, we are interested in transferring measures through slice-matching schemes as introduced in \cite{pitie2007automated}. We derive a generalized framework for these types of schemes and establish an interpretation as ``compatible maps'', which in turn allows for direct application to supervised learning tasks in the Wasserstein space. The main result of this paper is an a.s.\ convergence proof of a stochastic variant of the slice-matching schemes of \cite{pitie2007automated}, using stochastic gradient descent iterations in the Wasserstein space as suggested by \cite{backhoff2025stochastic}.
This convergence result contributes towards efforts in justifying the use of slice-matching schemes in data science applications \cite{bonneel2015sliced,rabin2010ShapeRetrieval,rabin2012wasserstein,meng2019large}. To this end, we show numerical experiments on image morphing.

\appendix

\section{{Key lemmas in the proof of \Cref{mainthm-Pj}}}\label{appendixProofs}

We consider a joint probability space for the sequence of random variables $\{P_k\}$, where $P_k\overset{\textrm{i.i.d.}}{\sim}u_n$ and  $u_n$ is the Haar probability measure on $O(n)$. The almost sure convergence result in \Cref{mainthm-Pj} is in terms of this joint probability space coupled with the product sigma-algebra and the associated product measure \cite{Saeki1996-uh}.

We show two lemmas, which are analogs of Proposition 3.2 and Lemma 3.3 in \cite{backhoff2025stochastic} adapted to our functional defined in \eqref{eq:LOjfunctional}. To quantify the behavior of ${H_j}(\sigma_k)$, given the first $k$ randomly chosen orthogonal matrices,  we introduce (in a similar fashion to \cite{backhoff2025stochastic}), the following (natural) filtration associated with the stochastic process $\{P_k\}_{k\in \mathbb{N}}$:
\begin{definition}\label{Pkdef}
    Let $P_k \overset{\textrm{i.i.d}}{\sim} u_n, k\geq 0$, where $u_n$ is the Haar probability measure on $O(n)$. Define $\mathcal{F}_0$ as the trivial the $\sigma$-algebra and $\mathcal{F}_{k+1}$ as the $\sigma$-algebra generated by $P_0,\ldots,P_k$.
\end{definition}

  \begin{lemma}\label{analog3.2} Let $ \sigma_{k+1}= \left((1-\gamma_k)id +\gamma_k\Pkjslice\right)_{\sharp}\sigma_k$ as defined in \eqref{sigmaiter-generalized}, where $P_k, \mathcal{F}_k$ are as in Definition \ref{Pkdef}. Then
      \begin{equation}
    \E[{H_j}(\sigma_{k+1})\mid \mathcal{F}_k]\leq (1+\gamma_k^2){H_j}(\sigma_k)-\gamma_k\|{H_j^\prime}(\sigma_k)\|^2_{\sigma_k}.
\end{equation} 
   \end{lemma}
   \begin{proof}
    Based on \Cref{lm:jsliceBound}, an analogous argument to Proposition 3.2 in \cite{backhoff2025stochastic} yields the following
       \begin{align*}
     &{H_j}(\sigma_{k+1})= \frac{1}{2}\int_{O(n)}\sum_{i=1}^j W_2^2(\sigma_{k+1}^{\theta_i}, \mu^{\theta_i})du_n(P)\nonumber\\
 \leq & \frac{1}{2}\int_{O(n)}\|\id-\jkPslice\|^2_{\sigma_k}du_n(P)+\frac{\gamma_k^2}{2} \int_{O(n)}\|\id-\Pkjslice\|^2_{\sigma_k}du_n(P)\nonumber\\
 &-\gamma_k\left<\int_{O(n)}(\id-\jkPslice)du_n(P),\id-\Pkjslice\right>_{\sigma_k}\nonumber\\
 = & \frac{1}{2}\int_{O(n)}\sum_{i=1}^jW_2^2(\sigma_k^{\theta_i}, \mu^{\theta_i})du_n(P)+\frac{\gamma_k^2}{2}\int_{O(n)}\sum_{i=1}^jW_2^2(\sigma_k^{\theta_i^k}, \mu^{\theta_i^k})du_n(P)\nonumber\\
 &-\gamma_k \left< {H_j}(\sigma_k), \id-\Pkjslice\right>_{\sigma_k},\label{B2lasteq}
   \end{align*}
    where $P = [\theta_1,\ldots,\theta_n]$ is a generic orthogonal matrix and $P_k = [\theta_1^k,\ldots,\theta_n^k ]$ is the orthogonal matrix at step $k$ of the iteration scheme \eqref{sigmaiter-generalized}. Note that the first two terms of the last inequality follows from essentially \Cref{lem:single-slice-distances} for the case where $\gamma_k =1$, and the last term follows from the the definition of the formal Frech\'et derivative in \eqref{FrecheDer-Pj}. Rewriting the above using the definition of ${H_j}$ and the fact that $u_n$ is a probability measure, we have
   \begin{equation*}
      {H_j}(\sigma_{k+1})\leq{H_j}(\sigma_k)+\frac{\gamma_k^2}{2}\sum_{i=1}^jW_2^2(\sigma_k^{\theta_i^k}, \mu^{\theta_i^k})-\gamma_k \left<{H_j}^{\prime}(\sigma_k), \id-\Pkjslice\right>_{\sigma_k}.\nonumber
   \end{equation*}
   Based on the above inequality, the final estimation for the conditional expectation  $ \mathbb{E}[{H_j}(\sigma_{k+1})|\mathcal{F}_k]$ {follows from taking conditional expectation with respect to $\mathcal{F}_k$ and the fact that $P_k$ is independently sampled from this $\sigma$-algebra,} which parallels the reasoning used in the last inequality in the proof of Proposition 3.2 in \cite{backhoff2025stochastic}.
\end{proof}

\renewcommand{\labelenumi}{\roman{enumi})} 
\begin{lemma}\label{analog3.3}
    Let $\{\rho_m\} \subseteq \W_{2,ac}(\R^n)$ such that $\rho_m \xrightarrow{W_2} \rho\in \W_{2,ac}(\R^n) $. Then as $m\rightarrow \infty$, we get
 \begin{enumerate}[(a)] 
        \item ${H_j}(\rho_m)\rightarrow{H_j}(\rho)$ and
        \item $\|{H^{\prime}_j}(\rho_m)\|_{\rho_m}\rightarrow \|{H^{\prime}_j}(\rho)\|_{\rho}$.
    \end{enumerate}
\end{lemma}

\begin{proof}
The arguments are inspired by the proof of Lemma 3.3 in \cite{backhoff2025stochastic}, with necessary adjustments made for the specific functional ${H_j}$.  For the convenience of the reader, we repeat the required set-up.

For part i), it suffices to show that $SW_2(\rho_m,\mu)\rightarrow SW_2(\rho,\mu)$ by the second equality in \Cref{rmk:Ljrelation}. Indeed, by triangle inequality and \Cref{lm:SWproperty}, we have 
$|SW_2(\rho_m,\mu)- SW_2(\rho,\mu)|\leq SW_2(\rho_m,\rho)\leq \frac{1}{\sqrt{n}}W_2(\rho_m,\rho) \rightarrow 0$ as $m\rightarrow \infty$.\\
\\
For part ii),  third equality in {\Cref{rmk:Ljrelation}} gives  that  ${H^{\prime}_j}(\eta) = j{H^{\prime}}(\eta)$ for $\eta \in \W_{2,ac}(\R^n)$ where 
\begin{equation}
   {H^{\prime}}(\eta) = \int_{S^{n-1}}\theta\left(\theta\cdot x - T^{\theta}_{\eta}\left(\theta\cdot x\right)\right)du(\theta),
\end{equation}
and $T^{\theta}_{\eta}$ denotes the optimal transport map from $\eta^{\theta}$
to $\mu^{\theta}$. Hence it suffices to show that $\|{H^{\prime}}(\rho_m)\|_{\rho_m}\rightarrow \|{H^{\prime}}(\rho)\|_{\rho}$ as $m\rightarrow\infty$. By the Skorokhod's representation theorem, there exists a common probability space $(\Omega, \mathcal{G}, \mathbb{P})$ for random vectors $\{X_m\}$ of laws $\{\rho_m\}$ and $X$ of law $\rho$ such that $X_m$ converges $\mathbb{P}$-a.s.\ to $X$. It follows that for every $\theta$,  $T^{\theta}_{\rho_m}(\theta\cdot X_m)\rightarrow T^{\theta}_{\rho}(\theta\cdot X)$ $\mathbb{P}$-a.s. by \cite[Theorem 3.4]{CUESTAALBERTOS199772}. Note that 
\begin{equation}\label{Lprimerhom}
    \|{H^{\prime}}(\rho_m)\|_{\rho_m}= \int_{\Omega}\left\|\int_{S^{n-1}}\theta\left(\theta\cdot X_m - T^{\theta}_{\rho_m}\left(\theta\cdot X_m\right)\right)du(\theta)\right \|^2d\mathbb{P}.
\end{equation}
By the triangle inequality, to obtain $\|{H^{\prime}}(\rho_m)\|_{\rho_m}\rightarrow \|{H^{\prime}}(\rho)\|_{\rho}$, it suffices to show that 
\begin{equation}\label{ineq:outertrisuff}
    \int_{\Omega}\left\|\int_{S^{n-1}}\theta\left[\theta\cdot X_m - T^{\theta}_{\rho_m}\left(\theta\cdot X_m\right)- \left(\theta\cdot X - T^{\theta}_{\rho}\left(\theta\cdot X\right)\right)\right]du(\theta)\right\|^2d\mathbb{P}\rightarrow 0. 
\end{equation}
By Jensen's inequalilty and the triangle inequality, it suffices to show the following two {properties}: 
\begin{align}
  &\theta\left(\theta\cdot X_m \right) \rightarrow    \theta\left(\theta\cdot X \right)  &\textrm{~in~} L^2(S^{n-1}\times \Omega),\label{ineq:Xm}\\
  &\theta T^{\theta}_{\rho_m}\left(\theta\cdot X_m\right) \rightarrow    \theta T^{\theta}_{\rho}\left(\theta\cdot X\right)   &\textrm{~in~} L^2(S^{n-1}\times \Omega). \label{ineq:TrhoXm}
\end{align}
Since $\|\theta\|\leq 1$ and $u(\theta)$ is the uniform probability measure on $S^{n-1}$, it is enough to show that
\begin{align}
  &\int_{\Omega}\|X_m-X\|^2d\mathbb{P}\rightarrow 0, \label{ineq:XmL2}\\
  &\int_{\Omega}\int_{S^{n-1}}\|T^{\theta}_{\rho_m}\left(\theta\cdot X_m\right)-  T^{\theta}_{\rho}\left(\theta\cdot X\right)\|^2 du(\theta)d\mathbb{P}\rightarrow 0 \label{ineq:TrhoXmL2}.
\end{align}
For {\eqref{ineq:XmL2}}, observe that $\sup\limits_{m} \int_{\|X_m\|\geq r}\|X_m(\omega)\|^2d\mathbb{P}(\omega) =\sup\limits_{m} \int_{\|x\|\geq r}\|x\|^2d\rho_m(x) \rightarrow 0 $ as $r\rightarrow\infty$ by \cite[Theorem 7.12, part (ii)]{villani2003topics}. Hence by the Vitali convergence theorem, \Cref{ineq:XmL2} holds. Similarly,  we have \Cref{ineq:TrhoXmL2}  by Fubini's theorem and observing that
\begin{align*}
 &\sup\limits_{m}\int_{S^{n-1}}  \int_{\left| T^{\theta}_{\rho_m}\left(\theta\cdot X_m\right)\right|\geq r} \left| T^{\theta}_{\rho_m}\left(\theta\cdot X_m\right)\right|^2d\mathbb{P}du(\theta)\\
 =&\sup\limits_{m}\int_{S^{n-1}} \int_{|t|\geq r} |t|^2d\mu^{\theta}(t)du(\theta)\\
 =& \sup\limits_{m}\int_{S^{n-1}} \int_{|\theta \cdot x |\geq r} |\theta\cdot x|^2d\mu(x)du(\theta)\\
 \leq& \sup\limits_{m}\int_{\|x\|\geq r}\|x\|^2d\mu(x) \rightarrow 0, \textrm{~ as~} r\rightarrow \infty,
\end{align*}
where the first two equalities follow from the change of variables formula and the fact that $\left(T^{\theta}_{\rho_m}\left(\theta\cdot X_m\right)\right)_\sharp \mathbb{P}= \mu^{\theta}$, and the last inequalilty uses the set inclusion: $\{x: |\theta\cdot x| \geq r\}\subset \{x: \|x\|\geq r\}$.

\end{proof}

\subsection{Auxiliary lemmas}

\begin{lemma}\label{lm:jsliceBound}
Let $P, P_k\in O(n)$, where $P = [\theta_1, \cdots, \theta_n]$ and $\sigma_{k+1}$ be defined in \eqref{sigmaiter-generalized} using $P_k$ and for fixed $1\leq j \leq n$. Then we get
    \begin{align*}
        \sum_{i=1}^j W_2^2(\sigma_{k+1}^{\theta_i},\mu^{\theta_i})&\leq \|\id-T^j_{\sigma_k,P}\|^2_{\sigma_k} +\gamma_k^2\|\id-\Pkjslice\|^2_{\sigma_k}\\
         &-2\gamma_k\left<\id-T^j_{\sigma_k,P}, \id-\Pkjslice\right>_{\sigma_k}.\nonumber
    \end{align*}
\end{lemma}
\begin{proof}
Let $\sigma_{k,\mu;P,j} = (T^j_{\sigma_k,P})_\sharp \sigma_k$. We first note that by \Cref{remark:single-slice} (4), $\mu^{\theta_i}= \sigma_{k,\mu;P,j}^{\theta_i}$ for $1\leq i \leq j$.
Observe that $\sigma_{k,\mu;P,j}^{\theta_i} = (\mathcal{P}_{\theta_i}\circ T^j_{\sigma_k,P})_{\sharp}\sigma_k$  and by  definition, $\sigma_{k+1}^{\theta_i}= \Big(\mathcal{P}_{\theta_i}\circ\left((1-\gamma_k)\id + \gamma_k \Pkjslice\right)\Big)_{\sharp}\sigma_k$. Then,
    by the Lipschitz continuity of the push-forward operation  \cite[Equation (2.1)]{Ambrosio2013}, for $1\leq i \leq j$ we get
    \begin{align}
        W_2^2(\sigma_{k+1}^{\theta_i},\mu^{\theta_i}) &\leq \left\|\mathcal{P}_{\theta_i}\circ \Big( T^j_{\sigma_k,P}-\left((1-\gamma_k)\id + \gamma_k \Pkjslice\right)\Big)\right\|^2_{\sigma_k}\nonumber\\
        & \leq \left\|\mathcal{P}_{\theta_i}\circ \Big((\id-T^j_{\sigma_k,P})-\gamma_k(\id-\Pkjslice)\Big)\right\|^2_{\sigma_k}.\nonumber
    \end{align}
    Note that for an orthonormal basis $\{\theta_1,\ldots,\theta_n\}$ of $\R^n$ and a function $F\in L^2(\sigma)$ we have $\|F\|^2_{\sigma} = \sum_{i=1}^n \|\mathcal{P}_{\theta_i}\circ F\|^2_{\sigma}$. Therefore, for an orthonormal set $\{\theta_1,\ldots,\theta_j\}$ with $j\leq n$ we get $ \sum_{i=1}^j \|\mathcal{P}_{\theta_i}\circ F\|^2_{\sigma}\leq \|F\|^2_{\sigma}$. This implies
\[\sum_{i=1}^j W_2^2(\sigma_{k+1}^{\theta_i},\mu^{\theta_i})  \leq \left\|(\id-T^j_{\sigma_k,P})-\gamma_k(\id-\Pkjslice)\right\|^2_{\sigma_k}.\]
      
The desired inequality hence follows from expanding the above inequality.
\end{proof}

\section{Propositions about the assumption (A1-i)}
{The following proposition gives a sufficient condition on which measures $\{\sigma_k\}$ produced in the proposed iterated scheme remains compactly supported (and absolutely continuous).  Consequently, they lie in a $W_2$-compact set, supporting Assumption $(A1-i)$; see \Cref{compact_W2} for more details. } 
\begin{proposition}\label{prop:A1_cond_1}
Let $\sigma_0\in \W_{2,ac}(\R^n)$. 
    Let $\tilde\Omega\subseteq \R^n$ be a compact set and $\mu\in \W_{2,ac}(\tilde\Omega)$. Then $\{\sigma_k\}_{k\geq 1}\subseteq \W_{2,ac}(\Omega)$ for some compact set $ \Omega \supseteq \tilde \Omega$, where $\sigma_k$'s are the defined as in the iterative scheme \Cref{sigmaiter-generalized}.
\end{proposition}
\begin{proof}
 Let $\tilde\Omega\subset B_r$, where $B_r$ is the closed ball of radius $r$ centered at the origin. Since $\mu^{\theta}$ has support on $[-r,r]$,  the essential range for the map $\Pkjslice$ is $B_R$, where $R=\sqrt{n}r$. The fact that $\sigma_k$'s are absolutely continuous follows from \cite[Lemma C.10]{limoos2024approximation}.
\end{proof}

{The previous proposition requires that the target $\mu$  be compactly supported. In the case of a non-compact case, the next two propositions and lemma show that a bounded moment condition on the measures provides an alternative sufficient condition for Assumption (A1-i).}
\begin{proposition}\label{prop:uniform3rdMoment}
 Let $\sigma_0,\mu\in \W_{2,ac}(\R^n)$ and $M:= \max\{M_3(\sigma_0), n^{1/2}M_3(\mu)\}<\infty$. For simplicity, let $\{\sigma_k\}_{k\geq 0}$ be defined as in \eqref{sigmaiter-generalized} for $j=n$. Then $\{M_3(\sigma_k)\}_{k\geq 0}\leq M$.  Here $M_3(\eta)$ denotes the $3$rd moment of a measure $\eta$.
\end{proposition}
\begin{proof}
 We will show that 
    \begin{align}\label{ineq_M3}
        M_3(\sigma_1)\leq \left(\left(1-\gamma_0\right )M_3(\sigma_0)^{1/3}+\gamma_0\left(n^{1/2}M_3(\mu)\right)^{1/3}\right)^3,
    \end{align}
    which immediately implies that $M_3(\sigma_1)\leq M$ and also  $M_3(\sigma_k)\leq M$ using induction. For simplicity denote $P_0 = [\theta_1, \cdots, \theta_n]\in O(n)$ and $f_{i}$ the optimal transport map between $\sigma^{\theta_i}$ to  $\mu^{\theta_i}$. Then $\sigma_1 = \Big(\sum_{i=1}^n\theta_if_i(\theta_i\cdot \ )\Big)_\sharp\sigma_0$. To see \eqref{ineq_M3},  
    \begin{align*}
        \int \|y\|^3 d\sigma_1(y)& = \int \|(1-\gamma_0)x + \gamma_0 \sum_{i=1}^nf_{i}(x\cdot \theta_i)\theta_i\|^3 d\sigma_0(x)\\
        &\leq \int \left((1-\gamma_0)\|x\|+ \gamma_0 \left(\sum_{i=1}^n f_{i}^2(x\cdot \theta_i)\right)^{1/2}\right)^3 d\sigma_0(x)\\
        &\leq \left((1-\gamma_0) \left(\int \|x\|^3 d\sigma_0(x)\right)^{1/3}d\sigma_0(x) + \gamma_0 \left(\int \left(\sum_{i=1}^nf_{i}^2(x\cdot \theta_i)\right)^{3/2}d\sigma_0(x)\right)^{1/3}\right)^3\\
        &\leq \left((1-\gamma_0)M_3^{1/3}(\sigma_0)+\gamma_0 \left(n^{1/2}M_3(\mu)\right)^{1/3}\right)^3\\
        &\leq (M^{1/3})^3 = M,
    \end{align*}
    where the first two inequalities follow from the triangle inequality and Minkowski inequality respectively, and the third inequality follows from {\Cref{lm:thirdMomentbd}}.
\end{proof}

\begin{lemma}\label{lm:thirdMomentbd}
 Assume that $\sigma_0,\mu\in \W_{2,ac}(\R^n)$ have finite 3rd moments. Then
    \begin{equation}
        \int \left(\sum_{i=1}^nf_{i}^2(x\cdot \theta_i)\right)^{3/2}d\sigma_0(x)\leq n^{1/2}M_3(\mu),
    \end{equation}
    where $\theta_i$'s and $f_i$'s are as in \Cref{prop:uniform3rdMoment}.
\end{lemma}
\begin{proof}
    Using the inequality that $\sum_{i=1}^n|a_i|^p \leq n^{p-1}\sum_{i=1}^n|a_i|^p$ for any $p\geq 1$, we have that 
    \begin{align*}
     \int \left(\sum_{i=1}^nf_{i}^2(x\cdot \theta_i)\right)^{3/2}d\sigma_0(x)&\leq n^{1/2} \int \sum_{i=1}^n |f_{i}(x\cdot \theta_i )|^3 d\sigma_0(x)\\
     &= n^{1/2}\sum_{i=1}^n \int_{\R} |f_{i}(t)|^3d\sigma^{\theta_i}(t)\\
     &= n^{1/2}\sum_{i=1}^n \int_{\R} |w|^3d\mu^{\theta_i}(t)\\
     &= n^{1/2}\sum_{i=1}^n\int |z\cdot \theta_i|^3d\mu(z)\\
     &\leq n^{1/2}\int \|z\|^3d\mu(z)=  n^{1/2}M_3(\mu).
    \end{align*}
\end{proof}

\begin{proposition}\label{prop:MomentSetCpt}
Let $K_{\sigma_0,\mu}:= \{\nu\in \W_{2,ac}(\R^n): \int_{\R^n}\|x\|^3d\nu(x)\leq M\}$. Then the closure of $K_{\sigma_0,\mu}$ is compact in $\W_{2}(\R^n)$.
\end{proposition}
\begin{proof}
 We first show that $K_{\sigma_0,\mu}$ is tight. Indeed, let $\nu\in K_{\sigma_0,\mu}$ and by Markov's inequality,
 \begin{align*}
    \int_{\|x\|\geq r} d\nu(x)&\leq \frac{\int \|x\|^3d\nu(x)}{r^3}\leq \frac{M}{r^3}.
 \end{align*}
  Let $\{\nu_j\}\subseteq K_{\sigma_0,\mu}$. Then by Prokhorov's theorem, $\nu_j \rightarrow \nu$ weakly for some $\nu\in \mathcal{P}(\R^n)$. Here we abuse notations and denote the convergent subsequence as $\{\nu_j\}$ as well. 
 Given $\epsilon>0$, let $r>0$ be such that $ \frac{M}{r}\leq \frac\epsilon2$. Let $J>0$ be such that for all $j\geq J$
 \begin{align*}
     \Biggl|\int_{B_r} \|x\|^2d\nu_j - \int_{B_r}\|x\|^2d\nu   \Biggr|< \frac{\epsilon}{2},
 \end{align*}
 where $B_r$ is the closed ball of radius $r$ centered at origin. Since $\int_{\R^n\setminus B_r} \|x\|^2d\nu_j \leq \int \frac{\|x\|^3}{r}d\nu_j \leq \frac{M}{r}\leq \frac{\epsilon}{2}$, we have that 
$\Biggl|\int \|x\|^2d\nu_j - \int_{B_r}\|x\|^2d\nu   \Biggr|<\epsilon$. Hence
  \begin{align*}
       M_2(\nu_j) \rightarrow M_2(\nu), ~\textrm{as~} j\rightarrow \infty.
 \end{align*}
 By \cite[Theorem 7.12]{villani2003topics}, $W_2(\nu_j,\nu)\rightarrow 0$ as $j\rightarrow \infty$.
\end{proof}

\section{Other technical details}

The following lemma is a special case of \cite[Lemma 5.2]{ma2023absolute}, which follows from \cite[Proposition 7.7, Remark 7.8]{santambrogio2015optimal} and references therein, see also \cite[Theorem 15.8, Theorem 15.9]{ambrosio2021lectures}.
\begin{lemma}(Special case of \cite[Lemma 5.2]{ma2023absolute})\label{lm:abs_continuity_energy}
Let $\leb$ be the Lebesgue measure on $\R^n$ and $G$ be a lower semi-continuous and convex function on $[0,\infty)$ such that 
\begin{align}\label{G_cond}
    G(t) \geq 0 \textrm{ and } \lim\limits_{t\rightarrow\infty} G(t)/t = \infty.
\end{align}
If a sequence $\{\rho_k := f_k\cdot \leb\}_k$ of absolutely continuous measures  converges weakly to a probability measure $\rho$ such that 
\begin{equation}\label{finite_G_energy}
    \liminf\limits_{k\rightarrow \infty}\int_{\R^n} G(f_k) d\leb \textrm{ is finite,}
\end{equation}
 then $\rho$ is also absolutely continuous with its density $f$ satisfying
\begin{equation}
    \int_{\R^n} G(f)d\leb  \leq    \liminf\limits_{k\rightarrow \infty} \int_{\R^n} G(f_k)d\leb.
\end{equation}
    
\end{lemma}

\begin{lemma}\label{sphereint_example}
Let $\{e_i\}_{i=1}^n$ be the standard basis for $\R^n$.
    \begin{equation*}
    \int_{S^{n-1}} (\theta\cdot e_i)\theta\,du(\theta) = [0,\cdots,0, w_i, 0\cdots, 0]^t,
\end{equation*}
where the $i$-th entry $w_i>0$.
\end{lemma}
\begin{proof}
Note that $\theta\cdot e_i = \theta(i)$.
    The natural symmetry inherent in $S^{n-1}$ allows us to readily observe that
    \begin{equation*}
    \int_{S^{n-1}}\theta(i)\theta(j)\,du(\theta) 
        \begin{cases}
           >0 \quad & i=j\\
           =0\quad & i\neq j
        \end{cases}. 
    \end{equation*}
\end{proof}

\begin{lemma}
   Let $\sigma\in \W_{2,ac}(\R^n)$. Then for any $\theta\in S^{n-1}$, $\sigma^{\theta}\in \W_{2,ac}(\R^n)$. 
\end{lemma}\label{absCont_proj}
 
\begin{proof}

The fact that $\sigma^{\theta}$ is absolutely continuous follows from the co-area formula, see also Box 2.4. in \cite[p. 82]{santambrogio2015optimal}
The fact that $\sigma^{\theta}$ has finite second moments follows from the change of variable formula
\begin{equation*}
    \int_{\R}|t|^2d\sigma^{\theta}(t)=\int_{\R^n}\|x\cdot\theta\|^2
    d\sigma(x)\leq \int_{\R^n}\|x^2\|d\sigma(x)<\infty. 
\end{equation*}
\end{proof}

\section{{Numerical Comparison with the Knothe-Rosenblatt transport}}\label{sec:KRcomparison}
{
A Python notebook for the following experiments is available in the accompanying GitHub Repository; see the {Data availability} section. }
{
\paragraph{Experiment D.1.}\label{Exp:rotGaussian}
 We compare the iterative matrix-slice-matching scheme $n=j$  with the Knothe--Rosenblatt (KR) transport on anisotropic Gaussian data in $\mathbb{R}^3$ ($n=3)$. Let the \emph{source} measure be $\sigma := \mathcal{N}(0,\Sigma)$ with $\Sigma=\mathrm{diag}(4,1,1/4)$, and the \emph{target} measure be $\mu := \mathcal{N}(0,Q\Sigma Q^\top)$ for a random orthogonal matrix $Q$. Each distribution is approximated with $N=200{,}000$ i.i.d.\ samples. For slice matching ($\gamma_k = \frac{1 + \log_{2}(k+2)}{\,k+2\,}, k\geq 0$), we apply the maps iteratively up to $20$ steps: after each iteration $k$ we evaluate the \emph{relative sliced-Wasserstein error} of the transported measure and record the \emph{cumulative runtime} up to iteration $k$. For KR, we report its (one-pass) error versus runtime. All curves are averaged over five trials.\\}
 \begin{figure}[h!]
    \centering
    \includegraphics[width=0.7\linewidth]{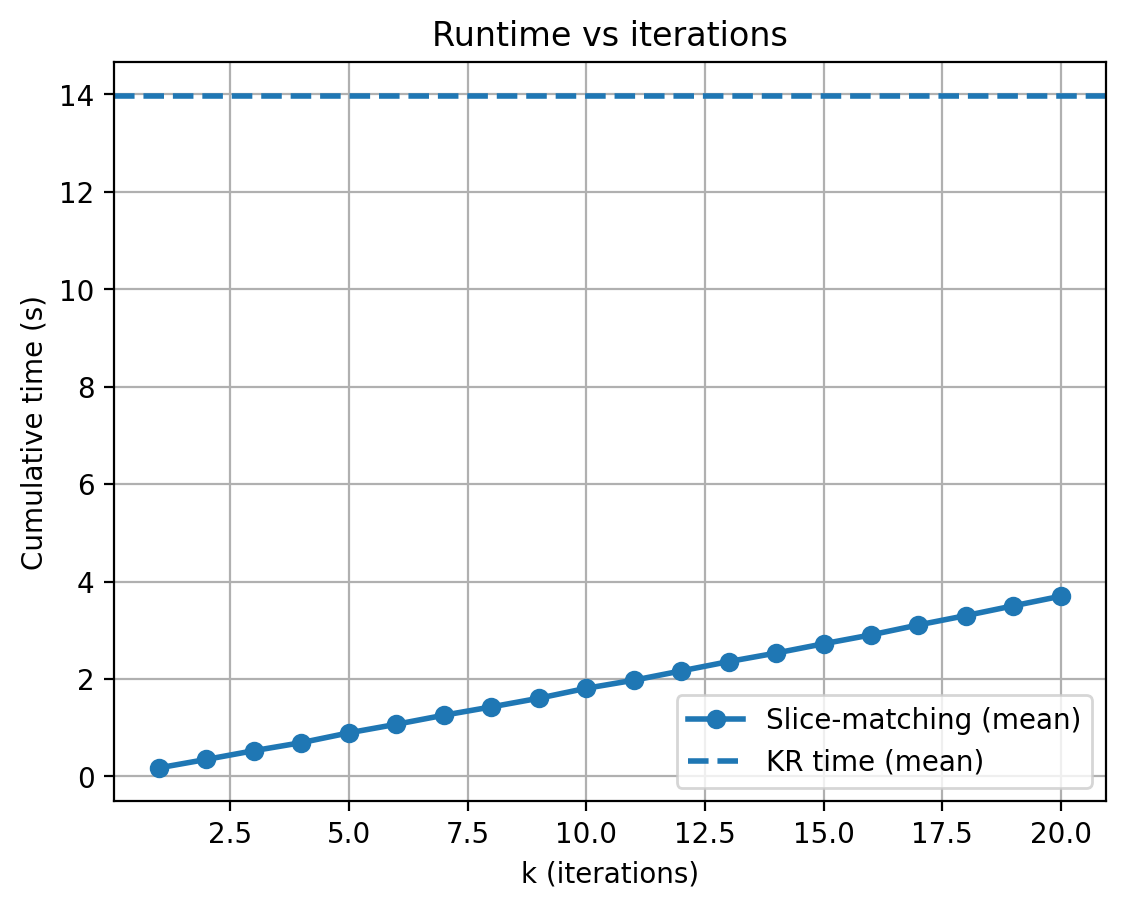}
    \caption{Cumulative runtime versus the number of iterations for Experiment D.1 (rotated Gaussians, $n=3)$.}
    \label{KRfig:n3rotatedGaus-iter}
\end{figure}

 { \paragraph{Discussion.} As expected, KR attains a small error in a single application; in this experiment, the mean relative error $\approx 1.49\times 10^{-3}$ and the mean runtime $\approx 14\,\mathrm{s}$. See \Cref{KRfig:n3rotatedGau-relError} for the error-runtime curves. Slice matching reaches a comparable error by the final iteration (mean relative error at $k=20$: $\approx 1.83\times 10^{-3}$) with substantially lower runtime (mean $\approx 4\,\mathrm{s}$), indicating a favorable error--time tradeoff. \Cref{KRfig:n3rotatedGaus-iter} further shows the mean cumulative runtime versus the iteration $k$ (averaged over five trials);  the KR curve is the average one-pass runtime shown as a constant reference.}

{
\paragraph{Experiment D.2.} We repeat the comparison on a \emph{shifted} Gaussian pair in dimension $n=20$. The source is $\sigma := \mathcal{N}(0,I_n)$ and the target is $\mu := \mathcal{N}(m,I_n)$ with $m=(50,0,\ldots,0)$. Each distribution is approximated with $N=5000$ i.i.d.\ samples. Slice matching ($\gamma_k = \frac{1 + \log_{2}(k+2)}{\,k+2\,}, k\geq 0$) is applied iteratively up to $k=20$; after each iteration we record the relative sliced–Wasserstein error of the transported measure and the cumulative runtime. For KR, we report its one-pass error versus runtime. All curves are averaged over five trials.\\} 

{\paragraph{Discussion.} Both KR and slice matching  attain near machine-precision accuracy in a single application. Note here $\gamma_0 = 1$, so in the matrix-slice-matching ($n=j$) case the target is recovered in one step, as illustrated in \Cref{rem:convergence-rate}.  The runtime of slice matching remains substantially smaller. See \Cref{KRfig:n20shiftedGau-relError} (error vs.\ runtime) and \Cref{KRfig:n20shiftedGau-iter} (cumulative runtime vs. iteration).}

\begin{figure}[h!]
    \centering
    \includegraphics[width=0.7\linewidth]{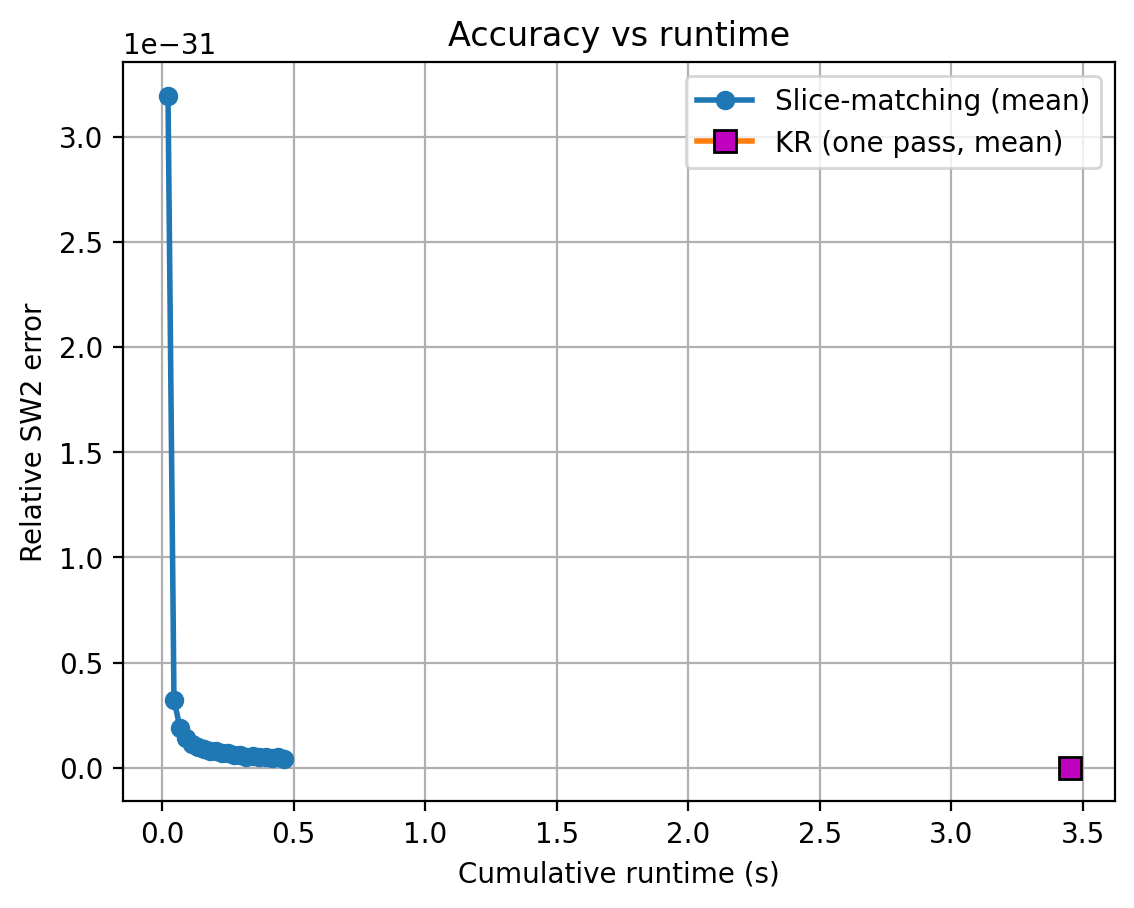}
    \caption{Comparison of measure transfer using slice-matching iterative scheme and Knothes-Rosenblatt transport (Experiment D.2). Source: $\mathcal{N}(0, I_n), n=20$. Target: $\mu := \mathcal{N}(m,I_n)$ with $m=(50,0,\ldots,0)$. Each distribution is approximated using  $N=5{,}000$ i.i.d. samples. The plot shows the relative sliced-Wasserstein error (averaged over 5 trials) versus  cumulative time. }
    \label{KRfig:n20shiftedGau-relError}
\end{figure}
\begin{figure}[h!]
    \centering
    \includegraphics[width=0.7\linewidth]{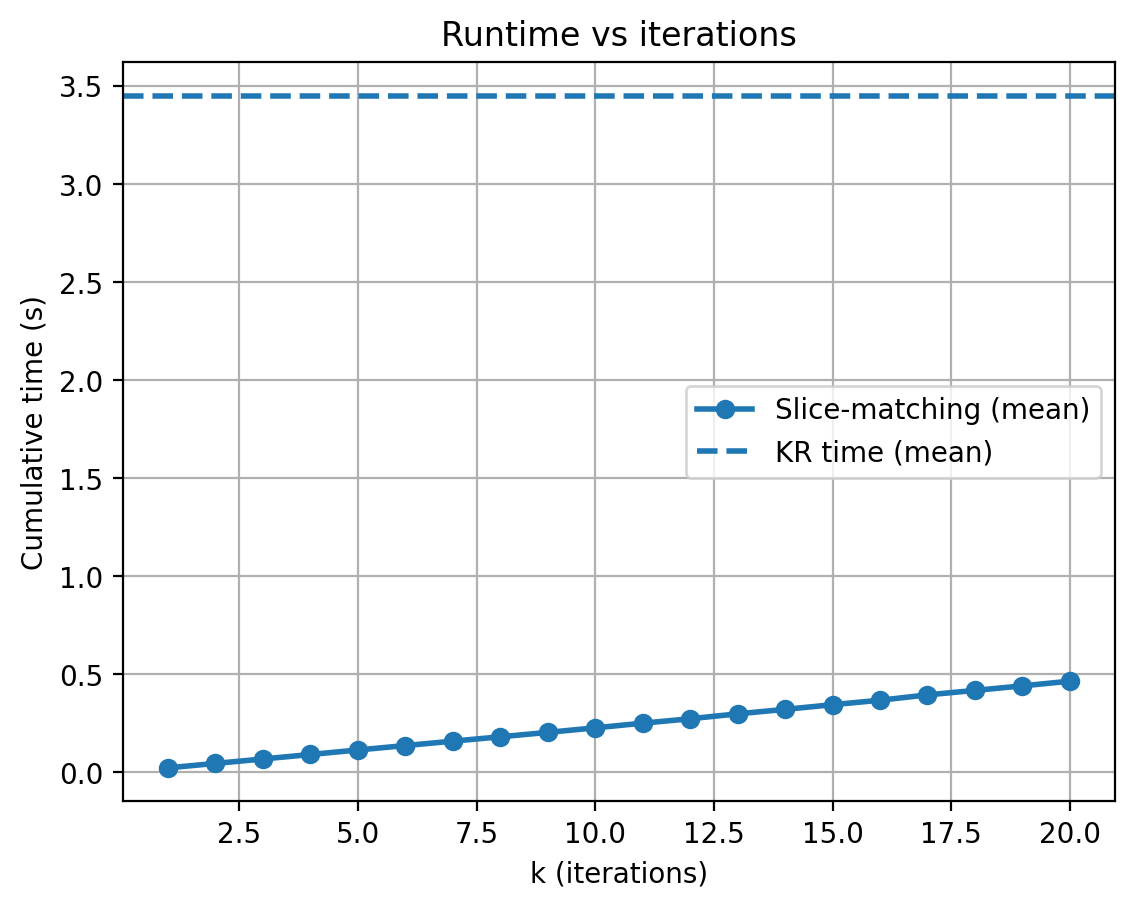}
    \caption{Cumulative runtime versus the number of iterations for Experiment D.2 (shifted Gaussians, $n=20)$.}
    \label{KRfig:n20shiftedGau-iter}
\end{figure}

{
\section{Additional numerical example}\label{sec:additional-numerics}

We expand the numerical section of \Cref{sec:experiments}. We show the digit merging experiment ($ 5\to 1$) of \Cref{sec:matrix-num-convergence}, but instead of using the matrix-slicing scheme ($j=2$), in every step, we only use a single slicing direction ($j=1$). We use the iterative scheme \eqref{sigmaiter-generalized} (note that $n=2$) and choice $\gamma_k = \frac{1+\log_2(k)}{k}$. 
In \Cref{fig:digit_morph_theta_log}, the error (in sliced Wasserstein distance) is shown as a function of the iteration variable $k$. Comparing to the matrix version ($j=n=2$) of \Cref{fig:digit_morph_matrix_log}, we see that the error decays slower in the single-slice case, which is to be expected.

\begin{figure*}[t!]
\centering
\includegraphics[scale=0.5]{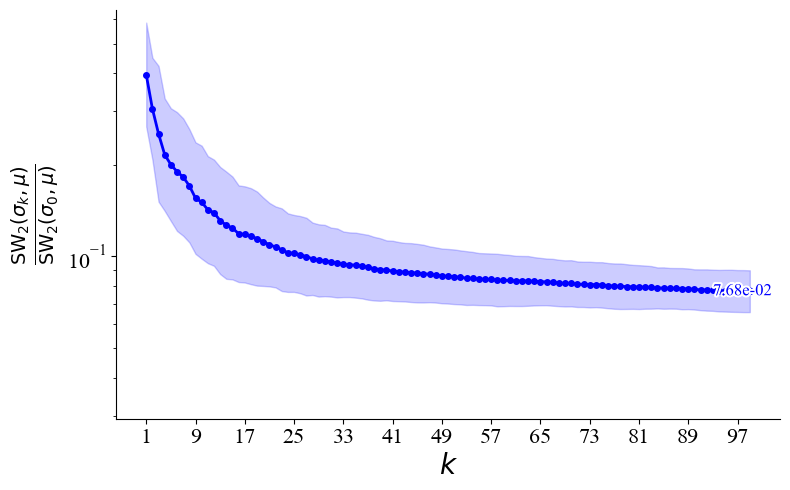}
\caption{Digit morphing experiment ($5 \to 1$) at $(84 \times 84)$ resolution, single-slice version ($j=1,n=2$) with 
log-decaying step size $\gamma_k=\frac{1+\log_2(k)}{k}$. 
Compared to the matrix version ($j=n=2$) of \Cref{fig:digit_morph_matrix_log}, we see that the error decays slower in the single-slice case.}
\label{fig:digit_morph_theta_log}
\end{figure*}

}

\section*{Acknowledgments}
Special thanks are extended to Soheil Kolouri for introducing SL to the convergence question related to the single-slice matching scheme. 
This work was partially supported by the National Science Foundation [DMS-2410140 to S.L. and C.M.; DMS-2306064 to C.M.]. C.M was also supported by a seed grant from the School of Data Science and Society at UNC and a faculty development award from UNC. 

\section*{Data availability}\label{DataAva}
The sample {Python code and notebooks} for this work are available at the \href{https://github.com/yongzhe-wang/slice_optimal_transport}{SliceMatchingTransport} Github repository.

\bibliographystyle{siam}
\bibliography{references}

@Article{mccann01,
    author="McCann, R.J.",
    title="Polar factorization of maps on {R}iemannian manifolds",
    journal="Geometric {\&} Functional Analysis GAFA",
    year="2001"
}

@article{cozzi2025long,
  title={Long-time asymptotics of the sliced-{W}asserstein flow},
  author={Cozzi, Giacomo and Santambrogio, Filippo},
  journal={SIAM Journal on Imaging Sciences},
  volume={18},
  number={1},
  pages={1--19},
  year={2025},
  publisher={SIAM}
}

@article{backhoff2025stochastic,
  title={Stochastic gradient descent for barycenters in {W}asserstein space},
  author={Backhoff, Julio and Fontbona, Joaquin and Rios, Gonzalo and Tobar, Felipe},
  journal={Journal of Applied Probability},
  volume={62},
  number={1},
  pages={15--43},
  year={2025},
  publisher={Cambridge University Press}
}

@article{baptista2023representation,
  title = {On the Representation and Learning of Monotone Triangular Transport Maps},
  author = {Baptista, Ricardo and Marzouk, Youssef and Zahm, Olivier},
  journal = {Foundations of Computational Mathematics},
  year = {2023},
  publisher = {Springer}
}

@book{williams1991probability,
  title={Probability with martingales},
  author={Williams, David},
  year={1991},
  publisher={Cambridge university press}
}

@article{bigot2018characterization,
  title={Characterization of barycenters in the Wasserstein space by averaging optimal transport maps},
  author={Bigot, J{\'e}r{\'e}mie and Klein, Thierry},
  journal={ESAIM: Probability and Statistics},
  volume={22},
  pages={35--57},
  year={2018},
  publisher={EDP Sciences}
}

@article{robbins1951stochastic,
  title={A stochastic approximation method},
  author={Robbins, Herbert and Monro, Sutton},
  journal={The {A}nnals of {M}athematical {S}tatistics},
  pages={400--407},
  year={1951},
  publisher={JSTOR}
}

@book{meckes_2019, place={Cambridge}, series={Cambridge Tracts in Mathematics}, title={The Random Matrix Theory of the Classical Compact Groups}, DOI={10.1017/9781108303453}, publisher={Cambridge University Press}, author={Meckes, Elizabeth S.}, year={2019}, collection={Cambridge Tracts in Mathematics}}

@article{pitie2007automated,
  title={Automated colour grading using colour distribution transfer},
  author={Piti{\'e}, Fran{\c{c}}ois and Kokaram, Anil C and Dahyot, Rozenn},
  journal={Computer Vision and Image Understanding},
  volume={107},
  number={1-2},
  pages={123--137},
  year={2007},
  publisher={Elsevier}
}

@inproceedings{frogner2015WassersteinLoss,
 author = {Frogner, Charlie and Zhang, Chiyuan and Mobahi, Hossein and Araya, Mauricio and Poggio, Tomaso A},
 booktitle = {Advances in Neural Information Processing Systems},
 editor = {C. Cortes and N. Lawrence and D. Lee and M. Sugiyama and R. Garnett},
 pages = {},
 publisher = {Curran Associates, Inc.},
 title = {Learning with a {W}asserstein Loss},
 volume = {28},
 year = {2015}
}

@article{limoos2024approximation,
  title={Approximation properties of slice-matching operators},
  author={Li, Shiying and Moosm{\"u}ller, Caroline},
  journal={Sampling Theory, Signal Processing, and Data Analysis},
  volume={22},
  number={1},
  pages={15},
  year={2024},
  publisher={Springer}
}

@book{ambrosio2008gradient,
  title={Gradient Flows: In Metric Spaces and in the Space of Probability Measures},
  author={Ambrosio, Luigi and Gigli, Nicola and Savare, Giuseppe},
  year={2008},
  publisher={Springer Science \& Business Media}
}

@article{benamou2000Fluid,
    author = {Benamou, Jean-David and Brenier, Yann},
    title = {A computational fluid mechanics solution to the {M}onge-{K}antorovich mass transfer problem},
    journal = {Numer.\ Math.},
    year = {2000},
    volume = {84},
    pages = {375--393},
    doi = {10.1007/s002110050002}
}

@article{santambrogio2017flows,
author = {Santambrogio, F.},
title = {\{Euclidean, metric, and {W}asserstein\} gradient flows: an overview},
journal = {Bull.\ Math.\ Sci.},
year = {2017},
volume = {7},
pages = {87--154}
}

@article{Jordan1998JKO,
author = {Jordan, Richard and Kinderlehrer, David and Otto, Felix},
title = {The Variational Formulation of the {F}okker--{P}lanck Equation},
journal = {SIAM Journal on Mathematical Analysis},
volume = {29},
number = {1},
pages = {1-17},
year = {1998},
doi = {10.1137/S0036141096303359},
}

@article{brenier1991,
  author  = {Brenier, Y.},
  title   = {Polar factorization and monotone rearrangement of vector-valued functions},
  journal = {Commun. Pure Appl. Math.},
  volume  = 44,
  number  = 4,
  pages   = {375--417},
  year    = 1991
}

@book{ambrosio2021lectures,
  title={Lectures on optimal transport},
  author={Ambrosio, Luigi and Bru{\'e}, Elia and Semola, Daniele and others},
  volume={130},
  year={2021},
  publisher={Springer}
}

@article{ma2023absolute,
  title={Absolute continuity of {W}asserstein barycenters on manifolds with a lower {R}icci curvature bound},
  author={Ma, Jianyu},
  journal={arXiv preprint arXiv:2310.13832},
  year={2023}
}

@article{santambrogio2015optimal,
  title={Optimal transport for applied mathematicians},
  author={Santambrogio, Filippo},
  journal={Birk{\"a}user, NY},
  volume={55},
  number={58-63},
  pages={94},
  year={2015},
  publisher={Springer}
}

@article{meng2019large,
  title={Large-scale optimal transport map estimation using projection pursuit},
  author={Meng, Cheng and Ke, Yuan and Zhang, Jingyi and Zhang, Mengrui and Zhong, Wenxuan and Ma, Ping},
  journal={Advances in Neural Information Processing Systems},
  volume={32},
  year={2019}
}

@article{zemel2019frechet,
  title={Fr{\'e}chet means and Procrustes analysis in Wasserstein space},
  author={Zemel, Yoav and Panaretos, Victor M},
  journal = {Bernoulli},
  volume = {25},
  number = {2},
  pages = {932--976},
  year={2019}
}

@phdthesis{bonnotte13thesis,
author = {Bonnotte, Nicolas},
year = {2013},
month = {12},
pages = {},
title = {Unidimensional and Evolution Methods for Optimal Transportation},
school = {Universit{\'e} Paris-Sud, Scuola Normale Superiore}
}

@article{bonneel2015sliced,
  title={Sliced and {R}adon {W}asserstein barycenters of measures},
  author={Bonneel, Nicolas and Rabin, Julien and Peyr{\'e}, Gabriel and Pfister, Hanspeter},
  journal={Journal of Mathematical Imaging and Vision},
  volume={51},
  pages={22--45},
  year={2015},
  publisher={Springer}
}

@article{bonet2022efficient,
  title={Efficient gradient flows in sliced-{W}asserstein space},
  author={Bonet, Cl{\'e}ment and Courty, Nicolas and Septier, Fran{\c{c}}ois and Drumetz, Lucas},
  journal={Transactions on Machine Learning Research},
  year={2022}
}

@inproceedings{chewi2020BuresWasserstein,
  title = {Gradient descent algorithms for {B}ures-{W}asserstein barycenters},
  author = {Chewi, Sinho and Maunu, Tyler and Rigollet, Philippe and Stromme, Austin J.},
  booktitle = {Proceedings of Thirty Third Conference on Learning Theory},
  pages = {1276--1304},
  year = {2020},
  editor = {Abernethy, Jacob and Agarwal, Shivani},
  volume = {125},
  series ={Proceedings of Machine Learning Research},
  publisher ={PMLR}
}

@inproceedings{mahey2023slicedGeodesics,
 author = {Mahey, Guillaume and Chapel, Laetitia and Gasso, Gilles and Bonet, Cl\'{e}ment and Courty, Nicolas},
 booktitle = {Advances in Neural Information Processing Systems},
 editor = {A. Oh and T. Naumann and A. Globerson and K. Saenko and M. Hardt and S. Levine},
 pages = {35350--35385},
 publisher = {Curran Associates, Inc.},
 title = {Fast Optimal Transport through Sliced Generalized {W}asserstein Geodesics},
 url = {https://proceedings.neurips.cc/paper_files/paper/2023/file/6f1346bac8b02f76a631400e2799b24b-Paper-Conference.pdf},
 volume = {36},
 year = {2023}
}

@INPROCEEDINGS{Bai2023partial,
  author={Bai, Yikun and Schmitzer, Bernhard and Thorpe, Matthew and Kolouri, Soheil},
  booktitle={2023 IEEE/CVF Conference on Computer Vision and Pattern Recognition (CVPR)}, 
  title={Sliced Optimal Partial Transport}, 
  year={2023},
  volume={},
  number={},
  pages={13681-13690},
  doi={10.1109/CVPR52729.2023.01315}}

@article{negrini2024applications,
  title={Applications of no-collision transportation maps in manifold learning},
  author={Negrini, Elisa and Nurbekyan, Levon},
  journal={SIAM Journal on Mathematics of Data Science},
  volume={6},
  number={1},
  pages={97--126},
  year={2024},
  publisher={SIAM}
}

@article{baptista2023ApproxFramework,
author = {Baptista, Ricardo and Hosseini, Bamdad and Kovachki, Nikola B. and Marzouk, Youssef M. and Sagiv, Amir},
title = {An Approximation Theory Framework for Measure-Transport Sampling Algorithms},
year = {2025},
journal = {Mathematics of Computation},
volume = {94},
number = {354},
pages = {1863 -- 1909},
doi = {10.1090/mcom/4013}
}

@INPROCEEDINGS{Kolouri2022ICASSP,
  author={Kolouri, Soheil and Nadjahi, Kimia and Shahrampour, Shahin and Şimşekli, Umut},
  booktitle={ICASSP 2022 - 2022 IEEE International Conference on Acoustics, Speech and Signal Processing (ICASSP)}, 
  title={Generalized Sliced Probability Metrics}, 
  year={2022},
  volume={},
  number={},
  pages={4513-4517},
  doi={10.1109/ICASSP43922.2022.9746016}}

@article{kolouri2019generalizedSlicesW,
author = {Kolouri, Soheil and Nadjahi, Kimia and Simsekli, Roland Badeau and Rohde, Gustavo},
title = {Generalized sliced {W}asserstein distances},
journal = {Advances in Neural Information Processing Systems},
volume = {32},
year = {2019}
}

@article{bonneel2019SPOT,
author = {Bonneel, Nicholas and Coeurjolly, David},
title = {SPOT: sliced partial optimal transport},
journal = {ACM Transactions on Graphics},
volume = {38},
number = {4},
pages = {1--13},
year = {2019}
}

@article{wang2013,
  Author  = {Wang, W. and Slep{\v{c}}ev, D. and Basu, S. and Ozolek, J.A. and Rohde, G.K.},
  Journal = {IJCV},
  Number  = {2},
  Pages   = {254--269},
  Title   = {A linear optimal transportation framework for quantifying and visualizing variations in sets of images},
  Volume  = {101},
  Year    = {2013}
}

@article{belkhatir2022texture,
title = {Wasserstein-based texture analysis in radiomic studies},
journal = {Computerized Medical Imaging and Graphics},
volume = {102},
pages = {102129},
year = {2022},
doi = {https://doi.org/10.1016/j.compmedimag.2022.102129},
author = {Zehor Belkhatir and Raúl San José Estépar and Allen R. Tannenbaum},
}

@inproceedings{sejourne2023unbalanced,
title={Unbalanced Optimal Transport meets Sliced-{W}asserstein},
author={Thibault Sejourne and Cl{\'e}ment Bonet and Kilian FATRAS and Kimia Nadjahi and Nicolas Courty},
booktitle={ICML Workshop on New Frontiers in Learning, Control, and Dynamical Systems},
year={2023},
}

@article{li2022neural,
title = {Sliced {W}asserstein Distance for Neural Style Transfer},
journal = {Computers \& Graphics},
volume = {102},
pages = {89-98},
year = {2022},
issn = {0097-8493},
doi = {10.1016/j.cag.2021.12.004},
author = {Jie Li and Dan Xu and Shaowen Yao},
}

@article{kolouri2017optimal,
  title={Optimal mass transport: Signal processing and machine-learning applications},
  author={Kolouri, Soheil and Park, Se Rim and Thorpe, Matthew and Slepcev, Dejan and Rohde, Gustavo K},
  journal={IEEE Signal Processing Magazine},
  volume={34},
  number={4},
  pages={43--59},
  year={2017},
  publisher={IEEE}
}

@inproceedings{solomon2014wasserstein,
  title={Wasserstein propagation for semi-supervised learning},
  author={Solomon, Justin and Rustamov, Raif and Guibas, Leonidas and Butscher, Adrian},
  booktitle={International Conference on Machine Learning},
  pages={306--314},
  year={2014}
}

@article{rubner2000earth,
  title={The earth mover's distance as a metric for image retrieval},
  author={Rubner, Yossi and Tomasi, Carlo and Guibas, Leonidas J},
  journal={Int J Comput Vis},
  volume={40},
  number={2},
  pages={99--121},
  year={2000},
  publisher={Springer}
}

@inbook{Ambrosio2013,
	address = {Berlin, Heidelberg},
	author = {Ambrosio, Luigi and Gigli, Nicola},
	pages = {1--155},
	publisher = {Springer Berlin Heidelberg},
	title = {A User's Guide to Optimal Transport},
	year = {2013}}

@article{kantorovich42,
    author = {Leonid Kantorovich},
    title = {On the transfer of masses},
    journal = {Doklady Akademii Nauk},
    volume = {37},
    number = {2},
    pages = {227-229},
    year = {1942}
}

@article{khurana2022supervised,
  title={Supervised learning of sheared distributions using linearized optimal transport},
  author={Khurana, Varun and Kannan, Harish and Cloninger, Alexander and Moosm{\"u}ller, Caroline},
  journal = {Sampling Theory, Signal Processing, and Data Analysis},
  volume = {21},
  number = {1},
  year={2023}
}

@article{moosmueller2020linear,
  title={Linear Optimal Transport Embedding: Provable {W}asserstein classification for certain rigid transformations and perturbations},
  author={Moosm{\"u}ller, Caroline and Cloninger, Alexander},
  journal={Information and Inference: A Journal of the IMA},
  volume={12},
  number={1},
  pages={363--389},
  year={2023}
}

@article{sinkhorn67,
    author = {Sinkhorn, R and Knopp, P},
    title = {Concerning nonnegative matrices and doubly stochastic matrices},
    journal = {Pacific J. Math.},
    volume = {21},
    number = {2},
    pages = {343-348},
    year = {1967}
}

@article{rosenblatt52arrangement,
 author = {Murray Rosenblatt},
 journal = {The Annals of Mathematical Statistics},
 number = {3},
 pages = {470--472},
 publisher = {Institute of Mathematical Statistics},
 title = {Remarks on a Multivariate Transformation},
 volume = {23},
 year = {1952}
}

@article{Knothe1957convex,
  title={Contributions to the theory of convex bodies},
  author={Herbert Knothe},
  journal={Michigan Mathematical Journal},
  year={1957},
  volume={4},
  pages={39-52}
}

@book{Villani1,
	Author = {Villani, C.},
	Publisher = {Springer-Verlag Berlin Heidelberg},
	Series = {Grundlehren der mathematischen Wissenschaften},
	Title = {Optimal Transport: Old and New},
	Volume = {338},
	Year = {2009}}

@book{villani2003topics,
  title={Topics in Optimal Transportation},
  author={Villani, C{\'e}dric},
  number={58},
  year={2003},
  publisher={American Mathematical Soc.}
}

@incollection{cuturi-2013,
    title = {Sinkhorn Distances: Lightspeed Computation of Optimal Transport},
    author = {Cuturi, M.},
    booktitle = {Advances in Neural Information Processing Systems 26},
    editor = {C. J. C. Burges and L. Bottou and M. Welling and Z. Ghahramani and K. Q. Weinberger},
    pages = {2292--2300},
    year = {2013},
    publisher = {Curran Associates, Inc.}
}

@article{nurbekyan2020nocollision,
    author = {Nurbekyan, L. and Iannantuono, A. and Oberman, A. M.},
    title = {No-Collision Transportation Maps},
    journal = {J Sci Comput},
    volume = {82},
    pages = {45},
    year = {2020}
}

@article{park2018cumulative,
  title={The cumulative distribution transform and linear pattern classification},
  author={Park, Se Rim and Kolouri, Soheil and Kundu, Shinjini and Rohde, Gustavo K},
  journal={Applied and computational harmonic analysis},
  volume={45},
  number={3},
  pages={616--641},
  year={2018},
  publisher={Elsevier}
}

@inproceedings{liutkus2019slicedFlows,
author = {Liutkus, Antoine and Simsekli, Umut and Majewski, Szymon and Durmus, Alain and St{\"o}ter, Fabian-Robert},
title = {Sliced-{W}asserstein
flows: Nonparametric generative modeling via optimal transport and diffusion},
year = {2019},
booktitle = {International Conference on
Machine Learning},
pages = { 4104--411}
}

@InProceedings{arjovsky2017wasserstein, 
  title = {{W}asserstein Generative Adversarial Networks}, 
  author = {Martin Arjovsky and Soumith Chintala and L{\'e}on Bottou},  
  pages = {214--223}, 
  year = {2017}, 
  editor = {Doina Precup and Yee Whye Teh}, 
  volume = {70}, 
  booktitle = {Proceedings of Machine Learning Research}, 
  publisher = {PMLR}
}

@article{aldroubi20,
  title = {Partitioning signal classes using transport transforms for data analysis and machine learning},
  journal = {Sampl. Theory Signal Process. Data Anal.},
  volume = {19},
  number = {6},
  author = {Akram Aldroubi and Shiying Li and Gustavo K. Rohde},
  year = {2021},
  doi = {10.1007/s43670-021-00009-z}
}

@article{CUESTAALBERTOS199772,
	author = {J.A. Cuesta-Albertos and C. Matr{\'a}n and A. Tuero-Diaz},
	journal = {Journal of Multivariate Analysis},
	number = {1},
	pages = {72-83},
	title = {Optimal Transportation Plans and Convergence in Distribution},
	volume = {60},
	year = {1997}}

@InProceedings{rabin2010ShapeRetrieval,
author="Rabin, Julien and Peyr{\'e}, Gabriel and Cohen, Laurent D.",
editor="Daniilidis, Kostas and Maragos, Petros and Paragios, Nikos",
title="Geodesic Shape Retrieval via Optimal Mass Transport",
booktitle="Computer Vision -- ECCV 2010",
year="2010",
publisher="Springer Berlin Heidelberg",
address="Berlin, Heidelberg",
pages="771--784",
}

@article{paulin2020,
    author = "Paulin, Loïs and Bonneel, Nicolas and Coeurjolly, David and Iehl, Jean-Claude and Webanck, Antoine and Desbrun, Mathieu and Ostromoukhov, Victor",
    title = "Sliced Optimal Transport Sampling",
    journal = "{ACM} Transactions on Graphics (Proceedings of SIGGRAPH)",
    year = "2020",
    volume = "39",
    number = "4",
    month = July
}

@inproceedings{
dai2021sliced,
title={Sliced Iterative Normalizing Flows},
author={Biwei Dai and Uros Seljak},
booktitle={ICML Workshop on Invertible Neural Networks, Normalizing Flows, and Explicit Likelihood Models},
year={2021},
url={https://openreview.net/forum?id=VmwEpdsvHZ9}
}

@article{Molchanov02_steepest,
	author = {Molchanov, Ilya and Zuyev, Sergei},
	journal = {Statistics and Computing},
	number = {2},
	pages = {115--123},
	title = {Steepest descent algorithms in a space of measures},
	volume = {12},
	year = {2002}}

@article{de1993new,
  title={New problems on minimizing movements},
  author={De Giorgi, Ennio},
  journal={Boundary value problems for PDE and applications},
  year={1993},
  publisher={Masson}
}

@ARTICLE{Saeki1996-uh,
  title     = "A proof of the existence of infinite product probability
               measures",
  author    = "Saeki, Sadahiro",
  journal   = "Am. Math. Mon.",
  publisher = "JSTOR",
  volume    =  103,
  number    =  8,
  pages     = "682",
  month     =  oct,
  year      =  1996
}

@inproceedings{rabin2012wasserstein,
  title={Wasserstein barycenter and its application to texture mixing},
  author={Rabin, Julien and Peyr{\'e}, Gabriel and Delon, Julie and Bernot, Marc},
  booktitle={Scale Space and Variational Methods in Computer Vision: Third International Conference, SSVM 2011, Ein-Gedi, Israel, May 29--June 2, 2011, Revised Selected Papers 3},
  pages={435--446},
  year={2012},
  organization={Springer}
}

\end{document}